\documentclass[amsfonts]{article}
\usepackage{amssymb}

\newtheorem{theorem}{Theorem}[section]

\newtheorem{example}[theorem]{Example}
\newtheorem{proposition}[theorem]{Proposition}
\newtheorem{lemma}[theorem]{Lemma}
\newtheorem{definition}[theorem]{Definition}
\newtheorem{remark}[theorem]{Remark}

\def\bR{\mathbb{R}}

\def\bN{\mathbb{N}}
\def\bD{\mathbb{D}}

\def\cB{\mathcal{B}}

\def\cF{\mathcal{F}}

\def\cH{\mathcal{H}}

\def\cP{\mathcal{P}}

\def\cS{\mathcal{S}}

\def\cX{\mathcal{X}}
\def\cY{\mathcal{Y}}

\begin{document}

\title{Stochastic Heat Equation
with Multiplicative Fractional-Colored Noise }

\author{Raluca M. Balan\thanks{Research supported by a grant from
the Natural Sciences and Engineering Research Council of Canada.}
$^{1}\quad$Ciprian A. Tudor\thanks{Research supported by the
Alexander von Humboldt Foundation.} $^{2}\vspace*{0.1in}$\\$^{1}$
Department of Mathematics and Statistics, University of Ottawa,
\\585, King Edward Avenue, Ottawa, ON, K1N 6N5, Canada.
\\rbalan@uottawa.ca\vspace*{0.1in}\\$^{2}$SAMOS/MATISSE,
Centre d'Economie de La Sorbonne,\\ Universit\'e de
Panth\'eon-Sorbonne Paris 1,\\90, rue de Tolbiac, 75634 Paris Cedex
13, France.\\tudor@univ-paris1.fr\vspace*{0.1in}}

\date{April 7, 2009}

\maketitle

\begin{abstract}
\noindent We consider the stochastic heat equation with
multiplicative noise $u_t=\frac{1}{2}\Delta u+ u \dot{W}$ in
$\bR_{+} \times \bR^d$, 
whose solution is interpreted in the mild sense. The noise $\dot W$
is fractional in time (with Hurst index $H \geq 1/2$), and colored
in space (with spatial covariance kernel $f$). When $H>1/2$, the
equation generalizes the It\^o-sense equation for $H=1/2$. We prove
that if $f$ is the Riesz kernel of order $\alpha$, or the Bessel
kernel of order $\alpha<d$, then the sufficient condition for the
existence of the solution is $d \leq 2+\alpha$ (if $H>1/2$),
respectively $d<2+\alpha$ (if $H=1/2$), whereas if $f$ is the heat
kernel or the Poisson kernel, then the equation has a solution for
any $d$. We give a representation of the $k$-th order moment of the
solution, in terms of an exponential moment of the ``convoluted
weighted'' intersection local time of $k$ independent
$d$-dimensional Brownian motions.
\end{abstract}




\section{Introduction}
The study of stochastic partial differential equations (s.p.d.e's)
driven by a Gaussian noise which is white in time and has a
non-trivial correlation structure in space (called ``color''),
constitutes now a classical line of research. These equations
represent an alternative to the standard s.p.d.e.'s driven by a
space-time white noise. A first step in this direction has been made
in \cite{DaFr98}, where the authors identify the necessary and
sufficient condition for the existence of the solution of the
stochastic wave equation (in spatial dimension $d=2$), in the space
of real-valued stochastic processes. The fundamental reference in
this area is Dalang's seminal article \cite{dalang99}, in which the
author gives the necessary and sufficient conditions under which
various s.p.d.e 's with a white-colored noise (e.g. the wave
equation, the damped heat equation, the heat equation) have a
process solution, in arbitrary spatial dimension. The methods used
in this article exploit the temporal martingale structure of the
noise, and cannot be applied when the noise is ``colored'' in time.
Other related references are: \cite{walsh86}, \cite{MiSa99},
\cite{PeZa00}, \cite{DaMu03} and \cite{DaSa05}.

Recently, there has been a growing interest in studying s.p.d.e.'s
driven by a Gaussian noise which has the covariance structure of the
fractional Browniam motion (fBm) in time, combined with a white (or
colored) spatial covariance structure. (Recall that an {\em fBm} is
a centered Gaussian process $(B_t)_{t \geq 0}$ with covariance
$E(B_tB_s)=R_{H}(t,s):= (t^{2H}+ s^{2H} -\vert t-s\vert ^{2H})/2$,
with $H \in (0,1)$. The Brownian motion is an fBm of index $H=1/2$.
We refer the reader to the expository article \cite{nualart03}, for
a comprehensive account on the fBm.) This interest comes from the
large number of applications of the fBm in practice. To list only a
few examples of the appearance of fractional noises in practical
situations, we mention \cite{KouSinney04} for biophysics,
\cite{BPS04} for financial time series, \cite{DMS03} for electrical
engineering, and \cite{CCL03} for physics.

In the present article, we consider the stochastic heat equation
with a multiplicative Gaussian noise, which is fractional (or white)
in time with Hurst index $H>1/2$ (respectively $H=1/2$), and has a
non-trivial spatial covariance structure given by a kernel $f$. As
in \cite{dalang99}, we assume that $f$ is the Fourier transform of a
tempered measure $\mu$. (Note that the particular case of a
spatially white noise arises when $f=\delta_0$.) More precisely, we
consider the following Cauchy problem:
\begin{eqnarray}
\label{heat} \frac{\partial u}{\partial t} &=& \frac{1}{2} \Delta u
+ u  \dot {W} , \quad t  >0, x\in \bR^{d} \\
\nonumber u_{0,x} &=& u_0(x), \quad x \in \bR^d,
\end{eqnarray}
where $u_{0} \in C_b(\bR^d)$ is non-random,
and $\dot W$ is a formal writing for the noise $W=\{W(h); h \in \cH
\cP\}$ (to be introduced rigourously in Section
\ref{prelim-section}).

Before discussing the multiplicative case, we recall briefly the
known results related to the existence of the solution of the
stochastic heat equation with additive noise:
\begin{eqnarray}
\label{heat-eq-HP} \frac{\partial u}{\partial t} &=& \frac{1}{2}
\Delta u + \dot {W} , \quad t  >0, x\in \bR^{d} \\
\nonumber u_{0,x} &=& 0, \quad x \in \bR^d,
\end{eqnarray}

When $H=1/2$ and $f=\delta_0$, equation (\ref{heat-eq-HP}) admits a
solution in the space of real-valued processes, if and only if
$d=1$. This phenomenon
can be explained intuitively by saying that, while the Laplacian
smooths, the white noise roughens (see also \cite{KFN08}). If the
spatial dimension $d$ is larger than $2$, then the roughness effect
of the white noise overcomes the smoothness influence of the
Laplacian.

What happens when the space-time white noise is replaced by a noise
which is fractional in time, but continues to be white in space?
This situation has been studied in several papers such as
\cite{DMD}, \cite{MaNu}, \cite{NuOuk}, \cite{TTV} and recently in
\cite{balan-tudor08}. In this case, a necessary and sufficient
condition for the existence of the solution of (\ref{heat-eq-HP}) is
$d<4H$, which allows us to consider the cases $d=1,2$ or $3$, for
suitable values of $H$. This can be interpreting by saying that for
$H>1/2$, the noise 
roughens a little bit less, and the smoothness influence of the
Laplacian overcomes the roughness of the noise. If the noise is
colored in space, the conditions for the existence of the solution
of (\ref{heat-eq-HP}) depend on the noise regularity in space. For
example, if $f$ the Riesz kernel of order $\alpha$, or the Bessel
kernel of order $\alpha$, the necessary and sufficient condition for
the existence of the solution of (\ref{heat-eq-HP}) is
$d<4H+\alpha$, whereas if $f$ is the heat kernel or the Poisson
kernel, the solution exists for any $d \geq 1$ and $H>1/2$ (see
\cite{balan-tudor08}, as well as Appendix B for a correction of the
result of \cite{balan-tudor08}).

Another explanation of this phenomenon is given in \cite{KFN08}, and
it is related to the local time of the stochastic processes
associated with the differential operator of the s.p.d.e. In the
particular case of the stochastic heat equation driven by a
space-time white noise, the solution exists only in dimension $d=1$
because this is the only case when the $d$-dimensional Brownian
motion has a local time.

We now return to the discussion of equation (\ref{heat}). This
equation has been studied recently in \cite{hu-nualart08}, when the
noise is fractional in time, and white in space. In this article, it
is proved that a {\em sufficient} condition for the existence of the
solution (in the space of square-integrable processes) is $d \leq
2$: if $d=1$, then equation (\ref{heat}) has a solution in any time
interval $[0,T]$, but if $d=2$, this equation has a solution only up
to a critical point $T_0$ (i.e. it has a solution in any interval
$[0,T]$, with $T<T_{0}$). It is not known if this condition is
necessary as well. There still is a connection with the local time,
in the sense that the second-order moment of the solution is equal
to the exponential moment of the ``weighted'' intersection local
time $L_{t}$ of two independent $d$-dimensional Brownian motion
$B^1$ and $B^2$, written formally as:
$$L_{t}:= H(2H-1)\int_{0}^{t}\int_{0}^{t} |r-s|^{2H-2}
\delta_{0}(B_{r}^{1}-B_{s}^{2})drds.$$

In the present article, we consider equation (\ref{heat}) driven by
the Gaussian noise introduced in \cite{balan-tudor08}. This noise is
fractional in time with Hurst index $H \geq 1/2$, and colored in
space, with covariance kernel $f$ chosen among the following: the
Riesz kernel, the Bessel kernel, the heat kernel, or the Poisson
kernel (see Examples \ref{riesz-kernel}-\ref{poisson-kernel}). The
case of the fractional kernel
$f(x)=\prod_{i=1}^{d}H_i(2H_i-1)|x_i|^{2H_i-2}$ with $1/2<H_i<1$ has
been examined in \cite{hu01}, using methods that rely on the product
form of $f$. These methods cannot be used in the present article,
since in our case (except the heat kernel), $f$ is not of product
type.  For the fractional kernel, it was proved in \cite{hu01} that
the sufficient condition for the existence of the solution
is $d<2/(2H-1)+\sum_{i=1}^d H_i$.

As in the case of equation (\ref{heat-eq-HP}) with additive noise,
we find that the existence of the solution depends on the roughness
of the noise. If $H>1/2$, and $f$ is the Riesz kernel of order
$\alpha$, or the Bessel kernel of order $\alpha<d$ (which are
``rough'' kernels), then a sufficient condition for the existence of
the solution is $d \leq 2+\alpha$: if $d<2+\alpha$ the solution
exists in any time interval $[0,T]$, whereas if $d=2+\alpha$, the
solution exists only up to a critical point $T_0$.
 If $f$ is the heat or the Poisson kernel (which are
``smooth'' kernels), the solution exists in any time interval, for
any $d \geq 1$ and $H \geq 1/2$. If $f$ is one of the ``rough''
kernels mentioned above, we prove that if the solution exists, then
$d<4H+\alpha$. This shows that for $H=1/2$, the necessary and
sufficient condition for the existence of solution is $d<2+\alpha$.
It remains an open problem to identify the necessary and sufficient
condition for the existence of the solution, in the case of $H>1/2$.

The existence of the solution is connected to the ``convoluted
weighted'' intersection local time $L_t$, written formally as:
$$L_t =H(2H-1) \int_{0}^{t} \int_{0}^{t} \int_{\bR^d}
|r-s|^{2H-2}\delta_{0}(B_{r}^{1}-B_{s}^{2}-y)f(y)dy dr ds.$$ More
precisely, the second-order moment of the solution can be expressed
as:
$$E[u_{t,x}^2]=E\left[u_0(x+B_t^1) u_0(x+B_t^2) \exp(L_t)\right].$$
As in \cite{hu-nualart08}, this expression can be extended to the
moments of order $k \geq 2$, using $k$ independent $d$-dimensional
Brownian motions.

This article is organized as follows. Section 2 contains some
preliminaries related to analysis on Wiener spaces. In Sections 3,
we discuss the existence of the solution. In Section 4, we examine
the relationship with the ``convoluted weighted'' intersection local
time.

\section{Preliminaries}
\label{prelim-section}

We begin by describing the kernel which gives the spatial covariance
of the noise. As in \cite{dalang99}, let $f$ be the Fourier
transform of a tempered distribution $\mu$ on $\bR^{d}$, i.e.
$$f(x)= \int_{\bR^{d}} e^{-i \xi \cdot x} \mu (d\xi ), \quad \forall
x \in \bR^d,$$ where $\xi \cdot x$ denotes the scalar product in
$\bR^d$. Let $\cP(\bR^d)$ be the completion of $\{1_{A}; A \in
\cB_b(\bR^d)\}$, where $\cB_{b}(\bR^d)$ denotes the class of bounded
Borel sets in $\bR^d$, with respect to the inner product
$$\langle 1_A, 1_B
\rangle_{\cP(\bR^d)}=\int_{A}\int_{B}f(x-y)dydx.$$

We consider some examples of kernel functions $f$. In what follows,
$\vert x\vert $ denotes the Euclidian norm of $x\in \mathbb{R}^d$.

\begin{example}
\label{riesz-kernel} {\rm  If $\mu(d\xi)=
 |\xi|^{-\alpha}d \xi$ for some $0<\alpha<d$, then $f$ is the Riesz kernel of order $\alpha$}:
$$f(x)=\gamma_{\alpha,d}|x|^{-d+\alpha},$$ where
$\gamma_{\alpha,d}=\Gamma((d-\alpha)/2)2^{-
\alpha}\pi^{-d/2}/\Gamma(\alpha/2)$.
\end{example}

\begin{example}
\label{bessel-kernel}{\rm If $\mu(d\xi)=(1+|\xi|^{2})^{-\alpha/2} d
\xi$ for some $\alpha>0$, then $f$ is the Bessel kernel of order
$\alpha$:
$$f(x)=\gamma_{\alpha}'\int_{0}^{\infty}w^{(\alpha-d)/2-1}
e^{-w}e^{-|x|^{2}/(4w)} dw,$$ where
$\gamma_{\alpha}'=(4\pi)^{\alpha/2}\Gamma(\alpha/2)$.
 In this case, $\cP(\bR^d)$ coincides
with $\cH^{-\alpha/2}(\bR^d)$, the fractional Sobolev space of order
$-\alpha/2$; see e.g. p.191, \cite{folland95}. }
\end{example}

\begin{example}
\label{heat-kernel} {\rm If $\mu(d\xi)=e^{-\pi^2 \alpha |\xi|^2/2} d
\xi$ for some $\alpha>0$, then $f$ is the heat kernel of order
$\alpha$:}
$$f(x)=(2 \pi \alpha)^{-d/2}
e^{-|x|^{2}/(2\alpha)}.$$
\end{example}

\begin{example}
\label{poisson-kernel} {\rm If $\mu(d\xi)=e^{-4\pi^2 \alpha |\xi|} d
\xi$ for some $\alpha>0$, then $f$ is the Poisson kernel of order
$\alpha$:}
$$f(x)=C_{d}\alpha
(|x|^{2}+\alpha^{2})^{-(d+1)/2},$$ where
$C_d=\pi^{-(d+1)/2}\Gamma((d+1)/2)$.
\end{example}

As in \cite{balan-tudor08}, if $H>1/2$, we let $\cH \cP$ be the
Hilbert space defined as the completion of $\{1_{[0,t] \times A}; t
\geq 0, A \in \cB_b(\bR^d) \}$ with respect to the inner product
\begin{equation}
\label{prosca} \langle 1_{[0,t] \times A}, 1_{[0,s] \times B}
\rangle_{\cH \cP}= \alpha_H \int_0^t \int_0^s \int_{A} \int_{B}
|u-v|^{2H-2}f(x-y)dy dxdv du,
\end{equation} where
$\alpha_H=H(2H-1)$. The space $\cH \cP$ is isomorphic to $\cH
\otimes \cP(\bR^d)$, where $\cH$ is the completion of $\{1_{[0,t]};
t \geq 0 \}$ with respect to the inner product
$$\langle 1_{[0,t]}, 1_{[0,s]} \rangle_{\cH}=\alpha_H
\int_0^t \int_0^s |u-v|^{2H-2}dv du.$$

If $H=1/2$, we let $\cH \cP$ be the completion of $\{1_{[0,t] \times
A}; t \geq 0, A \in \cB_b(\bR^d) \}$ with respect to the inner
product
$$\langle 1_{[0,t] \times A}, 1_{[0,s] \times B} \rangle_{\cH \cP}=
(t \wedge s) \int_{A} \int_{B} f(x-y)dy dx.$$ In this case, the
space $\cH \cP$ is isomorphic to $L^2(\bR_{+}) \otimes \cP(\bR^d)$.

We note that in both cases, the space $\cH \cP$ may contain
distributions.

\vspace{3mm}

Let $W=\{W(h); h \in \cH \cP\}$ be a zero-mean Gaussian process,
defined on a probability space $(\Omega, \cF,P)$, with covariance
$$E(W(h)W(g))=\langle h, g \rangle_{\cH \cP}.$$

The process $W$ introduce formally the noise perturbing the
stochastic heat equation. This noise is considered to be ``colored''
in space, with the color given by the kernel $f$. If $H>1/2$, the
noise is fractional in time, whereas if $H=1/2$ the noise is white
in time.

We now introduce the basic elements of analysis on Wiener spaces,
which are needed in the sequel. For a comprehensive account on this
subject, we refer the reader to \cite{nualart98} and
\cite{nualart06}.

We begin with a brief description of the multiple Wiener (or
Wiener-It\^o) 
 integral with respect to $W$. Let $\cF^{W}$ be the
$\sigma$-field generated by $\{W(h); h \in \cH \cP\}$, $H_n(x)$ be
the $n$-th order Hermite polynomial, and $\cH \cP_n$ be the closed
linear span of $\{H_n(W(h)); h \in \cH \cP\}$ in $L^2(\Omega,
\cF^{W},P)$. The space $\cH \cP_n$ is called the {\bf $n$-th Wiener
chaos} of $W$.

It is known that $L^2(\Omega, \cF^W,P)=\oplus_{n=0}^{\infty}\cH
\cP_n$, and hence every $F \in L^{2}(\Omega,\cF^W,P)$ admits the
following {\bf Wiener chaos expansion}:
\begin{equation}
\label{Wiener-chaos-0} F=\sum_{n=0}^{\infty}J_n(F),
\end{equation} where $J_n: L^2(\Omega,\cF^{W},P) \to \cH \cP_n$ is the
orthogonal projection. By convention, $\cH \cP_0=\bR$ and
$J_0(F)=E(F)$.

For each $n \geq 1$, and for each $h \in \cH \cP$ with $\|h \|_{\cH
\cP}=1$, we define
$$I_n(h^{\otimes n})=n! \ H_n(W(h)).$$
By polarization, we extend $I_n$ to elements of the form $h_1
\otimes \ldots \otimes h_n$ (see p. 230 of \cite{hu01}; e.g. $h_1
\otimes h_2= [(h_1+h_2)^{\otimes 2}-(h_1-h_2)^{\otimes 2}]/4$). By
linearity and continuity, we extend the definition of $I_n$ to the
space $\cH \cP^{\otimes n}$. (Note that if $\{e_i;i \geq 1\}$ is a
CONS in $\cH \cP$, then $\{e_{i_1} \otimes \ldots \otimes e_{i_n};
i_j \geq 1\}$ is a CONS in $\cH \cP^{\otimes n}$.) For any $h \in
\cH \cP^{\otimes n}$, we say that
$$I_n(h):=\int_{(\bR_{+} \times \bR^{d})^n}h(t_1,x_1,\ldots,
t_n,x_n)dW_{t_1,x_1} \ldots dW_{t_n, x_n}$$ is the {\bf multiple
Wiener integral} of $h$ with respect to $W$. We have
$$E(I_n(h)I_n(g))=n! \ \langle \tilde h, \tilde g \rangle_{\cH \cP^{\otimes
n}}, \quad \forall h,g \in \cH \cP^{\otimes n},$$ where $\tilde
h(t_1,x_1, \ldots,t_n,x_n)=(n!)^{-1}\sum_{\sigma \in S_n}
h(t_{\sigma(1)}, x_{\sigma(1)}, \ldots,
t_{\sigma(n)},x_{\sigma(n)})$ is the symmetrization of $h$ 
with respect to the $n$ variables $(t_{1}, x_{1}), \ldots ,
(t_{n},x_{n})$, and $S_n$ is the set of all permutations of
$\{1,\ldots,n\}$. By convention, we set $I_0(x)=x$.

The map $I_n: \cH \cP^{\otimes n} \to \cH \cP_n$ is surjective.
Moreover, for any $F_n \in \cH \cP_n$, there exists a unique $f_n
\in \cH \cP^{\otimes n}$ symmetric, such that $I_n(f_n)=F_n$. Using
(\ref{Wiener-chaos-0}), we conclude that any $F \in
L^2(\Omega,\cF^{W},P)$ can be written as:
\begin{equation}
\label{Wiener-chaos} F=\sum_{n=0}^{\infty}I_n(f_n),
\end{equation}
where $f_0=E(F)$ and $f_n \in \cH \cP^{\otimes n}$ is symmetric and
uniquely determined by $F$. We have:
$$E|F|^2=\sum_{n=0}^{\infty}
E|I_n(f_n)|^2=\sum_{n=0}^{\infty}n! \ \|f_n\|_{\cH \cP^{\otimes
n}}^2.$$

\vspace{3mm}

We now introduce the stochastic integral with respect to $W$. Let
$u=\{u_{t,x}; (t,x) \in \bR_{+} \times \bR^d\}$ be an
$\cF^W$-measurable square-integrable process. By
(\ref{Wiener-chaos}), for any $(t,x)\in \bR_{+} \times \bR^d$, we
have
\begin{equation}
\label{Wiener-chaos-u-tx}
u_{t,x}=E(u_{t,x})+\sum_{n=1}^{\infty}I_n(f_n(\cdot,t,x)),
\end{equation} where $f_n(\cdot, t,x) \in \cH \cP^{\otimes n}$ is
symmetric and uniquely determined by $u_{t,x}$. For each $n \geq 1$,
let $\tilde f_n$ be the symmetrization of $f_n$ with respect to all
$n+1$ variables. Let $\tilde f_0=E(u)$.
We say that $u$ is {\bf integrable with respect to $W$} if $\tilde
f_n \in \cH \cP^{\otimes (n+1)}$ for every $n \geq 0$, and
$\sum_{n=0}^{\infty}I_{n+1}(\tilde f_n)$ converges in $L^2(\Omega)$.
In this case, we define the stochastic integral
$$\delta(u):=\int_0^{\infty}u_s \delta W_s=\sum_{n=0}^{\infty}I_{n+1}(\tilde f_n)$$
Note that: $$E|\delta(u)|^2=\sum_{n=0}^{\infty}(n+1)! \ \|\tilde
f_n\|_{\cH \cP^{\otimes (n+1)}}^{2}.$$

The following alternative characterization of the operator $\delta$
is needed in the present article. Let $\cS=\{F=f(W(h_1),\ldots,
W(h_n)); f \in C_{b}^{\infty}(\bR^n), h_i \in \cH \cP, n \geq 1\}$
be the space of all ``smooth cylindrical'' random variables,  where
$C_b^{\infty}(\bR^d)$ denotes the class of all bounded infinitely
differentiable functions on $\bR^n$, whose partial derivatives are
also bounded. The Malliavin derivative of an element
$F=f(W(h_1),\ldots, W(h_n))\in \cS$, with respect to $W$, is defined
by:
$$DF:=\sum_{i=1}^{n}\frac{\partial f}{\partial x_i}(W(h_1),\ldots,
W(h_n))h_i.$$ Note that $DF \in L^2(\Omega; \cH \cP)$; by abuse of
notation, we write $DF=\{D_{t,x} F;(t,x) \in [0,T] \times \bR^d\}$
even if $D_{t,x} F$ is not a function in $(t,x)$. We endow $\cS$
with the norm $\|F\|_{\bD^{1,2}}^{2}:=E|F|^2+E\|D F \|_{\cH
\cP}^{2}$, we let $\bD^{1,2}$ be the completion of $\cS$ with
respect to this norm. The operator $D$ can be extended to
$\bD^{1,2}$. Then $\delta: \mbox{Dom} \ \delta \subset L^{2}(\Omega;
\cH \cP) \to L^2(\Omega)$ is the adjoint
of the operator $D$, and is uniquely defined by the following 
duality relationship: $u \in \mbox{Dom} \ \delta$ if and only if
\begin{equation}
\label{duality} E(F \delta(u))=E\langle DF,u \rangle_{\cH \cP},
\quad \forall F \in \bD^{1,2}.
 \end{equation}
Note that $u \in \mbox{Dom} \ \delta$ if and only if $u$ is
integrable with respect to $W$. (In the literature, $\delta$ is
called the Skorohod integral with respect to $W$.)


\section{Existence of the solution}
\label{existence-solution-section}

In this section, we give conditions for the existence of the
solution of equation (\ref{heat}).

Let $p_t(x)$ be the heat kernel on $\bR^d$, i.e.
$$p_{t}(x)=\frac{1}{(2\pi t)^{d/2}}\exp\left(-\frac{|x|^2}{2t} \right),
\quad t>0, x \in \bR^d.$$ For any bounded Borel function $\varphi:
\bR^d \to \bR$, let $p_t \varphi
(x)=\int_{\bR^d}p_{t}(x-y)\varphi(y)dy$.

\vspace{3mm}

For each $t >0$, let $\cF_{t}$ be the $\sigma$-field generated by
$\{W(1_{[0,s] \times A}); s \in [0,t], A \in \cB_b(\bR^d)\}$.

The solution of equation (\ref{heat}) is interpreted in the mild (or
evolution) sense, using the stochastic integral introduced above.
More precisely, we have the following definition.

\begin{definition}
An $(\cF_t)_t$-adapted square-integrable process $u=\{u_{t,x}; (t,x)
\in \bR_{+} \times \bR^{d}\}$ is a {\bf solution to} (\ref{heat}) if
for any $(t,x) \in \bR_{+} \times \bR^d$, the process
$$\{Y_{s,y}^{t,x}=1_{[0,t]}(s) p_{t-s}(x-y) u_{s,y} ; \ (s,y) \in \bR_{+} \times
\bR^d\}$$ is integrable with respect to $W$, and $$u_{t,x}=p_t
u_0(x)+\int_{0}^{\infty}\int_{\bR^d}Y_{s,y}^{t,x} \delta W_{s,y}.$$
\end{definition}

\noindent By (\ref{duality}), the above definition is equivalent to
saying that for any $(t,x) \in \bR_{+} \times \bR^{d}$, $u_{t,x} \in
L^{2}(\Omega)$, $u_{t,x}$ is $\cF_t$-measurable, and
\begin{equation}
\label{def-sol-u} E(u_{t,x}F)=E(F)p_t u_0(x)+E \langle Y^{t,x}, DF
\rangle_{\cH \cP}, \quad \forall F \in \bD^{1,2}.
\end{equation}

 The next result establishes the existence of the solution
$u=\{u_{t,x}; (t,x) \in \bR_{+} \times \bR^{d}\}$, as a collection
of random variables in $L^{2}(\Omega)$. As in \cite{hu-nualart08}
(see also \cite{buckdahn-nualart94}, \cite{leon-sanmartin07},
\cite{nourdin-tudor06}, \cite{nualart-rozovskii97},
\cite{nualart-zakai89}, \cite{tudor04}), one can find a closed
formula for the kernels $f_n(\cdot,t,x)$ appearing in the Wiener
chaos expansion (\ref{Wiener-chaos-u-tx}) of $u_{t,x}$.

\begin{proposition}
\label{existence} In order that equation (\ref{heat}) possesses a
solution it is necessary and sufficient that for any $ (t,x) \in
\bR_{+} \times \bR^d$, we have
\begin{equation}
\label{cond-fn-existence} \sum_{n=0}^{\infty}n! \
\|f_n(\cdot,t,x)\|_{\cH \cP^{\otimes n}}^2<\infty,
\end{equation}
where
$$f_{n}(t_{1},x_{1},\ldots,t_{n},x_{n},t,x) = \frac{1}{n!}
\prod _{j=1}^{n}p_{t_{\rho(j+1)}-t_{\rho(j)}}(x_{\rho(j+1)}
-x_{\rho(j)})p_{t_{\rho(1)}}u_0(x_{\rho(1)}),$$ $\rho$ denotes the
permutation of $\{1,2, \ldots, n\}$ such that
$t_{\rho(1)}<t_{\rho(2)}<\ldots <t_{\rho(n)}$, $t_{\rho(n+1)}=t$ and
$x_{\rho (n+1)}=x$. In this case, the solution $u$ is unique in
$L^2(\Omega)$, admits the Wiener chaos decomposition
(\ref{Wiener-chaos-u-tx}), and
\begin{equation}
\label{moment-2-utx} E|u_{t,x}|^2=\sum_{n=0}^{\infty}n! \
\|f_n(\cdot,t,x)\|_{\cH \cP^{\otimes n}}^2.
\end{equation}
\end{proposition}

We begin to examine condition (\ref{cond-fn-existence}). Note that
\begin{equation}
\label{estimate-alpha-n-prime} \alpha_n(t,x):=(n!)^2 \
\|f_n(\cdot,t,x)\|_{\cH \cP^{\otimes n}}^2 \leq \|u_0\|_{\infty}^2
\alpha_n(t),
\end{equation}
with equality if $u_0=1$, and hence
$$E|u_{t,x}|^2=\sum_{n=0}^{\infty}\frac{1}{n!}\alpha_n(t,x)\leq
\|u_0\|_{\infty}^2 \sum_{n=0}^{\infty}\frac{1}{n!}\alpha_n(t),$$
where
\begin{equation}
\label{def-alpha-n} \alpha_n(t) = \left\{
\begin{array}{ll} \alpha_H^n \int_{[0,t]^{2n}} \prod_{j=1}^{n}|s_j-t_j|^{2H-2}
\psi^{*(n)}({\bf s}, {\bf t}) d{\bf s} d {\bf t}
& \mbox{if $H>1/2$} \\
\\
\int_{[0,t]^{n}} \psi^{*(n)}({\bf s}, {\bf s}) d{\bf s} & \mbox{if
$H=1/2$}
\end{array} \right.
\end{equation}
and
\begin{eqnarray}
\psi^{*(n)}({\bf s}, {\bf t}) & := & \int_{\bR^{2nd}}
\prod_{j=1}^{n}f(x_j-y_j)\prod_{j=1}^{n}p_{t_{\rho(j+1)}-t_{\rho(j)}}(x_{\rho(j+1)}-x_{\rho(j)})
\nonumber \\
& & \prod_{j=1}^{n}p_{s_{\sigma(j+1)}-s_{\sigma(j)}}
(y_{\sigma(j+1)}-y_{\sigma(j)}) d {\bf x} d {\bf y}.\label{psin}
\end{eqnarray}
In the above integrals, we denoted ${\bf s}=(s_1,\ldots, s_n), {\bf
t}=(t_1,\ldots, t_n), {\bf x}=(x_1,\ldots,x_n), \linebreak {\bf
y}=(y_1,\ldots, y_n)$, and we chose the permutations $\rho$ and
$\sigma$ of $\{1,\ldots,n\}$ such that
\begin{equation}
\label{permut-rho-sigma} 0<t_{\rho(1)}<t_{\rho(2)}< \ldots <
t_{\rho(n)} \quad \mbox{and} \quad 0<s_{\sigma(1)}<s_{\sigma(2)}<
\ldots < s_{\sigma(n)},
\end{equation} with
$t_{\rho(n+1)}=s_{\sigma(n+1)}=t$ and
$x_{\rho(n+1)}=y_{\sigma(n+1)}=x$.

Note that $$\psi^{*(n)}({\bf s},{\bf t}) = \langle g_{\bf s}^{(n)},
g_{\bf t}^{(n)} \rangle_{\cP(\bR^d)^{\otimes n}}, \quad \forall {\bf
t},{\bf s} \in [0,t]^n,$$ where
\begin{eqnarray*}
g_{\bf t}^{(n)}(x_1, \ldots x_n) & = &
\prod_{j=1}^{n}p_{t_{\rho(j+1)}-
t_{\rho}(j)}(x_{\rho(j+1)}-x_{\rho(j)}) \\
 g_{\bf s}^{(n)}(y_1, \ldots
y_n) & = & \prod_{j=1}^{n} p_{s_{\sigma(j+1)}-
s_{\sigma}(j)}(y_{\sigma(j+1)}-y_{\sigma(j)}),
\end{eqnarray*}
and the permutations $\rho$ and $\sigma$ are chosen such that
(\ref{permut-rho-sigma}) holds.

\vspace{2mm}

As in \cite{balan-tudor08}, for any $y,z \in \bR^d$ and $u,v>0$, we
denote
$$J_{f}(u,v,y,z):=\int_{\bR^d} \int_{\bR^d}p_{u}(x-y)p_{v}(x'-z)f(x-x')dxdx'.$$

\begin{lemma}
\label{estimate-J} (i) If $f$ is the Riesz kernel of order $\alpha$,
or  the Bessel kernel of order $\alpha<d$, then
$$J_{f}(u,v,y,z) \leq D_{\alpha,d} (u+v)^{-(d-\alpha)/2}, \quad \forall y,z \in
\bR^d,$$ where $D_{\alpha,d}$ is a positive constant depending on
$\alpha$ and $d$.

(ii) If $f$ is the heat kernel of order $\alpha$, or the Poisson
kernel of order $\alpha$, then
$$J_{f}(u,v,y,z) \leq C_{\alpha,d}, \quad \forall y,z \in
\bR^d.$$
\end{lemma}

\noindent {\bf Proof:} Note that
$J_{f}(u,v,y,z)=E[f(y-z+\sqrt{u}Y-\sqrt{v}Z)]$, where $Y$ and $Z$
are independent $d$-dimensional standard normal random vectors. We
use the following inequality: (see (3.24) of \cite{balan-tudor08})
\begin{equation}
\label{chi-sq-ineq} E[e^{-|y-z+\sqrt{2u}Y-\sqrt{2v}Z|^2 /(4w)}] \leq
\left(1+ \frac{u+v}{w}\right)^{-d/2}.
\end{equation}

(i) In the case of the Riesz kernel, this inequality has been shown
in the proof of Theorem 3.13 of \cite{balan-tudor08}. Suppose now
that $f$ is the Bessel kernel of order $\alpha<d$. Using
(\ref{chi-sq-ineq}),
\begin{eqnarray*} J_{f}(u,v,y,z) &=&
\gamma_{\alpha}' \int_0^{\infty}
w^{(\alpha-d)/2-1}e^{-w}E[e^{-|y-z+\sqrt{u}Y-\sqrt{v}Z|^2/(4w)}]dw
\\
&\leq & \gamma_{\alpha}' \int_0^{\infty} w^{(\alpha-d)/2-1}e^{-w}
\left(1+\frac{u+v}{2w}\right)^{-d/2}dw \\
&=& \gamma_{\alpha}' \int_0^{\infty} w^{\alpha/2-1} e^{-w} \left(
w+\frac{u+v}{2}\right)^{-d/2}dw=\gamma_{\alpha}'I_{\alpha,d}\left(\frac{u+v}{2}\right) \\
\end{eqnarray*}
where $I_{\alpha,d}(x):=\int_{0}^{\infty}w^{\alpha/2-1} e^{-w} (
w+x)^{-d/2}dw$. The result follows, since
\begin{eqnarray*}
I_{\alpha,d}(x) & \leq & x^{-d/2} \int_0^x w^{\alpha/2-1}
e^{-w}dw+\int_{x}^{\infty}w^{\alpha/2-1} e^{-w} ( w+x)^{-d/2}dw \\
&=& x^{-d/2} x^{\alpha/2} \int_0^1 y^{\alpha/2-1}e^{-xy}dy +
x^{-(d-\alpha)/2}\int_1^{\infty}y^{\alpha/2-1}e^{-xy} (y+1)^{-d/2}dy
\\
&\leq & x^{-(d-\alpha)/2} \int_0^1 y^{\alpha/2-1} dy +
x^{-(d-\alpha)/2}\int_1^{\infty}y^{-(d-\alpha)/2-1}dy=K_{\alpha,d}x^{-(d-\alpha)/2},
\end{eqnarray*}
where $K_{\alpha,d}=2/\alpha+2/(d-\alpha)=2d/[\alpha(d-\alpha)]$,
and we used the fact that $\alpha<d$.

(ii) If $f$ is the heat kernel, using (\ref{chi-sq-ineq}), we
obtain:
\begin{eqnarray*}
J_{f}(u,v,y,z) &=& (2\pi
\alpha)^{-d/2}E[e^{-|y-z+\sqrt{u}Y-\sqrt{v}Z|^2/(2\alpha)}] \\
& \leq & (2\pi \alpha)^{-d/2} \left(1+
\frac{u+v}{\alpha}\right)^{-d/2} \\
&=& (2\pi)^{-d/2} (\alpha+u+v)^{-d/2} \leq (2\pi\alpha)^{-d/2}.
\end{eqnarray*}
If $f$ is the Poisson kernel, we have:
\begin{eqnarray*}
J_{f}(u,v,y,z) &=& C_{d}\alpha
E[(|y-z+\sqrt{u}Y-\sqrt{v}Z|^2+\alpha^2)^{-(d+1)/2}] \\
& \leq & C_d \alpha (\alpha^2)^{-(d+1)/2}=C_d\alpha^{-d}.
\end{eqnarray*}
$\Box$

\begin{lemma}
\label{estimate-psi} (i) If $f$ is the Riesz kernel of order
$\alpha$, or the Bessel kernel of order $\alpha <d$, then for any
${\bf s},{\bf t} \in [0,t]^n$,
$$ \psi^{*(n)}({\bf s}, {\bf t}) \leq \left(D_{\alpha,d}
2^{-(d-\alpha)/2} \right)^{n} [ \beta({\bf s})\beta({\bf
t})]^{-(d-\alpha)/4},$$ where $\beta({\bf s}):=\prod_{j=1}^{n}
(s_{\sigma(j+1)}- s_{\sigma(j)})$, $\beta({\bf t}): =
\prod_{j=1}^{n} (t_{\rho(j+1)}- t_{\rho(j)})$, and the permutations
$\rho$ and $\sigma$ are chosen such that (\ref{permut-rho-sigma})
holds.

(ii) If $f$ is the heat kernel of order $\alpha$, or the Poisson
kernel of order $\alpha$, then for any ${\bf s},{\bf t} \in
[0,t]^n$,
$$\psi^{*(n)}({\bf s}, {\bf t}) \leq C_{\alpha,d}^n,$$
where $C_{\alpha,d}$ is a constant depending on $\alpha$ and $d$.
\end{lemma}

\noindent {\bf Proof:} By the Cauchy-Schwartz inequality,
$$\psi^{*(n)}({\bf s}, {\bf t}) \leq
\psi^{*(n)}({\bf s}, {\bf s})^{1/2}\psi^{*(n)}({\bf t}, {\bf
t})^{1/2}.$$

\noindent To find an upper bound for $\psi^{*(n)}({\bf s}, {\bf
s})$, we use Lemma \ref{estimate-J} to estimate the following
integrals:
\begin{eqnarray*}
I_j&:=&\int_{\bR^d}\int_{\bR^d} p_{u_j}
(x_{\sigma(j+1)}-x_{\sigma(j)})
p_{u_j}(y_{\sigma(j+1)}-y_{\sigma(j)})
f(x_{\sigma(j)}-y_{\sigma(j)})dx_{\sigma(j)} dy_{\sigma(j)}\\
&=& J_{f}(u_j, u_j, x_{\sigma(j+1)},y_{\sigma(j+1)}), \quad
\mbox{with} \ u_j=s_{\sigma(j+1)}-s_{\sigma(j)}, \ j=1, \ldots,n.
\end{eqnarray*}

(i) In this case, $I_j \leq D_{\alpha,d}
[2(s_{\sigma(j+1)}-s_{\sigma(j)})]^{-(d-\alpha)/2}$ and
$$\psi^{*(n)}({\bf s}, {\bf s}) \leq D_{\alpha,d}^n 2^{-n(d-\alpha)/2}
\left[\prod_{j=1}^{n} (s_{\sigma(j+1)}- s_{\sigma(j)})
\right]^{-(d-\alpha)/2}.$$

(ii) In this case, $I_j \leq C_{\alpha,d}$ and $\psi^{*(n)}({\bf s},
{\bf s}) \leq C_{\alpha,d}^n$. $\Box$

\vspace{3mm}

If $H>1/2$, it was proved in \cite{MMV01} that there exists
$\beta_{H}>0$ such that
$$\alpha_{H}
\int_{0}^{\infty} \int_{0}^{\infty} \varphi
(s)\varphi(t)|t-s|^{2H-2}ds dt \leq \beta_H^2 \left(
\int_{0}^{\infty} |\varphi(t)|^{1/H}dt \right)^{2H},$$ for any
$\varphi \in L^{1/H}(\bR_{+})$. Hence,
\begin{equation}
\label{ineq-norms-H-L} \alpha_{H}^n \int_{\bR_{+}^{2n}} \varphi
({\bf s})\varphi({\bf t})\prod_{j=1}^{n}|t_j-s_j|^{2H-2}d {\bf s} d
{\bf t} \leq \beta_H^{2n} \left( \int_{\bR_+^n} |\varphi({\bf
t})|^{1/H}d{\bf t} \right)^{2H},
\end{equation}
for any $\varphi \in L^{1/H}(\bR_{+}^n)$. If $H=1/2$, we let
$\beta_H=1$.

\vspace{2mm}

We need the following auxiliary result.

\begin{lemma}
\label{calcul-In} Let
$$I_{n}(t,h):= \int_{T_n} [(t-s_n)(s_n-s_{n-1}) \ldots
(s_2-s_1)]^{h} d{\bf s},$$ where $T_n=\{{\bf s}=(s_1, \ldots,s_n);
0<s_1<s_2 < \ldots <s_n<t\}$. Then $I_n(t,h)<\infty$ if and only if
$1+h>0$. In this case,
$$I_n(t,h)=\frac{\Gamma(1+h)^{n+1}}{\Gamma(n(1+h)+1)}
t^{n(1+h)}.$$
\end{lemma}

\noindent {\bf Proof:} First note that $\int_{0}^{s_{2}}
(s_{2}-s_{1} )^{h} ds_{1}= s_{2}^{h+1}/(h+1)$, and then
\begin{eqnarray*}
\int_{0}^{s_{3}}(s_{3} -s_{2}) ^{h}s_{2} ^{h+1}ds_{2} &= &s_{3} ^{2h+2 } \beta \left( h+2, h+1\right)\\
&=& s_{3}^{2(h+1)} \beta \left(  (h+1)+1, h+1)\right)
\end{eqnarray*}
where $\beta(a,b)=\int_0^1x^{a-1}(1-x)^{b-1}dx$ is the beta
function, and we used the change of variables $s_{2}/s_{3}=z$. In
this way, the integral $I_{n}(t,h)$ becomes
\begin{eqnarray*}
I_{n}(t,h)
&=& \frac{1}{h+1} \beta \left( (h+1)+1, h+1\right) \ldots \beta \left( (n-2)(h+1) +1, h+1\right)\\
 &&  \int_{0}^{t}
s_{n} ^{(n-1)(h+1)}(t-s_{n})^{h}ds_{n} \\
&=& t^{n(h+1) } \beta \left( (h+1)+1, h+1\right) \ldots \beta \left(
(n-1) (h+1) +1, h+1\right).
\end{eqnarray*}
Using the fact that $\beta (a,b)= \Gamma (a) \Gamma (b) /\Gamma
(a+b)$ for $a,b>0$ and $\Gamma (z+1) = z\Gamma (z) $ for any $z>0$,
we obtain the desired conclusion. $\Box$

Using Lemma \ref{estimate-psi}, Lemma \ref{calcul-In} and
(\ref{ineq-norms-H-L}), we obtain the following estimate for
$\alpha_n(t)$.

\begin{proposition}
\label{estimate-alpha-n} Suppose that $H \geq 1/2$ and let
$\alpha_n(t)$ be given by (\ref{def-alpha-n}).

(i) If $f$ is the Riesz kernel of order $\alpha$, or the Bessel
kernel of order $\alpha <d$, and
\begin{equation}
\label{cond-d-alphaf-1} H>\frac{d-\alpha}{4} \ ,
\end{equation}
then
$$\alpha_n(t)  \leq C_{H,d,\alpha}^* D(t)^{n} \
(n!)^{(d-\alpha)/2}, \quad \mbox{for any} \ t>0, n \geq 1,$$ where
$C_{H,d,\alpha}^*>0$ is a constant depending on $H,d,\alpha$, and
$$D(t)=D_{\alpha,d} 2^{-(d-\alpha)/2} \beta_H^2 \
\Gamma\left(1-\frac{d-\alpha}{4H} \right)^{2H}
\left(1-\frac{d-\alpha}{4H}
\right)^{-[2H-(d-\alpha)/2]}t^{2H-(d-\alpha)/2}.$$

(ii) If $f$ is the heat kernel of order $\alpha$, or the Poisson
kernel of order $\alpha$, then
$$\alpha_n(t) \leq C(t)^n, \quad \mbox{for any} \ t>0
\ \mbox{and} \ n \geq 1,$$ where $C(t)=C_{\alpha,d} t^{2H}$.
\end{proposition}

\begin{remark}
{\rm Proposition \ref{estimate-alpha-n} shows that $\sum_n
\alpha_n(t)/n!$ grows exponentially in time in some cases, and
faster than exponentially in other cases.}
\end{remark}

\noindent {\bf Proof of Proposition \ref{estimate-alpha-n}:} We only
give the proof in the case $H>1/2$, the case $H=1/2$ being similar.
We use the definition (\ref{def-alpha-n}) of $\alpha_n(t)$.

(i) Let $h=-(d-\alpha)/(4H)$. By Lemma \ref{estimate-psi}.(i) and
(\ref{ineq-norms-H-L}), we obtain:
\begin{eqnarray*}
\alpha_n(t) & \leq & \left(D_{\alpha,d} 2^{-(d-\alpha)/2}
\right)^{n} \alpha_H^n \int_{([0,t]^{2})^{n}} \prod_{j=1}^{n}
|t_j-s_j|^{2H-2} [ \beta({\bf s})\beta({\bf t})]^{-(d-\alpha)/4}
d{\bf s}
d{\bf t} \\
& \leq & \left(D_{\alpha,d} 2^{-(d-\alpha_f)/2} \right)^{n}
\beta_H^{2n} \left( \int_{[0,t]^{n}} \beta({\bf
s})^{-(d-\alpha)/(4H)} d{\bf s}
\right)^{2H} \\
&=& \left(D_{\alpha,d} 2^{-(d-\alpha)/2} \beta_H^2 \right)^{n}
(n!)^{2H} I_n(t,h)^{2H}.
\end{eqnarray*}

\noindent Using Lemma \ref{calcul-In}, we obtain:
\begin{eqnarray*}
\alpha_n(t) & \leq &  \Gamma(1+h)^{2H}\left(D_{\alpha,d}
2^{-(d-\alpha)/2} \ \beta_H^2 \right)^{n} \  (n!)^{2H} \left\{
\frac{\Gamma(1+h )^{n}}{\Gamma(n(1+h)+1)}
t^{n(1+h)}  \right\}^{2H}\\
&=& \Gamma(1+h)^{2H} \left\{D_{\alpha,d} 2^{-(d-\alpha)/2} \
\beta_H^2 \Gamma\left(1+h \right)^{2H}  t^{2H(1+h)} \right\}^{n} \
\left(\frac{n!}{\Gamma(n(1+h)+1)}\right)^{2H}.
\end{eqnarray*}
The result follows by using (3.19) of \cite{hu-nualart08}.

(ii) By Lemma \ref{estimate-psi}.(ii),
$$\alpha_n(t)  \leq C_{\alpha,d}^n \alpha_H^n
\int_{[0,t]^{2n}} \prod_{j=1}^{n}|s_j-t_j|^{2H-2} d{\bf s} d{\bf t}
= C_{\alpha,d}^n \ t^{2Hn}=C(t)^n.$$ $\Box$

\vspace{3mm}

Using Proposition \ref{estimate-alpha-n}, we examine the existence
of the solution of equation (\ref{heat}). The next result is an
extension of Proposition 4.3 of \cite{hu-nualart08} to the case of a
colored noise $W$.

\begin{proposition}
\label{our-existence} (i) Let $f$ be the Riesz kernel of order
$\alpha$, or the Bessel kernel of order $\alpha <d$. Suppose that
either
\begin{equation}
\label{cond-d-alphaf-2} H>1/2 \quad \mbox{and} \quad d \leq
2+\alpha,
\end{equation}
or
\begin{equation}
\label{cond-d-alphaf-3} H=1/2 \quad \mbox{and} \quad d < 2+\alpha.
\end{equation}

\noindent Then (\ref{heat}) has a unique solution in $[0,T] \times
\bR^{d}$, provided that $T<T_0$ where
\begin{equation}
\label{def-of-T0} T_0= \left\{
\begin{array}{ll} \left\{\left(1-\frac{1}{2H}\right)D_{\alpha,d} 2^{-1} \beta_H^2 \Gamma
\left(1-\frac{1}{2H} \right)^{2H} \right\}^{-1/(2H-1)}
& \mbox{if $d=2+\alpha$} \\
\infty & \mbox{if $d<2+\alpha$}
\end{array} \right.
\end{equation}

(ii) Let $H \geq 1/2$, and $f$ be the heat kernel of order $\alpha$,
or the Poisson kernel of order $\alpha$. Then (\ref{heat}) has a
unique solution in $\bR_{+} \times \bR^d$.
\end{proposition}

\begin{remark}
{\rm Either one of conditions (\ref{cond-d-alphaf-2}) or
(\ref{cond-d-alphaf-3}) is stronger that (\ref{cond-d-alphaf-1}).}
\end{remark}

\begin{remark}
{\rm Proposition \ref{our-existence} shows that in the case $H=1/2$,
the dimension $d=2+\alpha$ cannot be attained. }
\end{remark}

\noindent {\bf Proof of Proposition \ref{our-existence}:} We apply
Proposition \ref{existence}, using Proposition
\ref{estimate-alpha-n}.

(i) We have:
$$\sum_{n=0}^{\infty}n! \ \|f_n(\cdot,t,x)\|_{\cH \cP^{\otimes n}}^2
\leq  \|u_0\|_{\infty}^2 \sum_{n=0}^{\infty} \frac{1}{n!}
\alpha_n(t) \leq  \|u_0\|_{\infty}^2 C_{H,\alpha,d}^*
\sum_{n=0}^{\infty} \frac{D(t)^{n}}{(n!)^{1-(d-\alpha)/2}}.$$

If $d-\alpha=2$, then the last sum is finite if $D(t)<1$, which is
equivalent to saying that $t<T_0$. If $d-\alpha<2$, then the last
sum is finite for any $t>0$, by Stirling's formula and D'Alembert
criterion.

(ii) We have:
$$\sum_{n=0}^{\infty}n! \ \|f_n(\cdot,t,x)\|_{\cH \cP^{\otimes n}}^2
\leq  \|u_0\|_{\infty}^2 \sum_{n=0}^{\infty} \frac{1}{n!}
\alpha_n(t) \leq \|u_0\|_{\infty}^2 \sum_{n=0}^{\infty}
\frac{C(t)^n}{n!} <\infty.$$ $\Box$

The next result shows that  (\ref{cond-d-alphaf-1}) is a necessary
condition for the existence of the solution.

\begin{proposition}
\label{neces-RB} Suppose that $H \geq 1/2$ and $f$ is either the
Riesz kernel or order $\alpha$, or the Bessel kernel of order
$\alpha$. If equation (\ref{heat}) with $u_0=1$ has a solution in
$\bR_{+} \times \bR^d$, then (\ref{cond-d-alphaf-1}) holds.
\end{proposition}

\noindent {\bf Proof:} Note that $E\left| u_{t,x} \right| ^{2}=
\sum_{n=0}^{\infty} \alpha_{n}(t)/n!<\infty$ implies that
$\alpha_1(t) <\infty$, which in turn implies (\ref{cond-d-alphaf-1})
(see Appendix A). $\Box$

\begin{remark}
\label{nec-remark} {\rm Proposition \ref{our-existence} and
Proposition \ref{neces-RB} show that, if $H=1/2$ and $f$ is the
Riesz kernel of order $\alpha$, or the Bessel kernel of order
$\alpha<d$, then the condition $d<\alpha+2$ is necessary and
sufficient for the existence of the solution of (\ref{heat}). It
remains an open problem to see if this condition is necessary, when
$H>1/2$. To resolve this issue, one needs to develop a full analysis
of the range of $\alpha_n(t)$, which would include the
identification of suitable lower bounds. Such analysis will be the
subject of future investigations. }
\end{remark}

\begin{remark}
{\rm The case $H<\frac{1}{2}$ also constitutes an interesting line
of investigation, which will be pursued in a subsequent article. We
mention that in this case even the stochastic heat equation with
linear additive fractional-colored noise has not been solved. The
technical difficulties that appear here are related to the structure
of the space ${\cal{HP}}$ and the lack of the expression of the
scalar product in this space as (\ref{prosca}). Indeed, when
$H<\frac{1}{2}$, assuming that the noise $W(t,x)$ is defined for
$t\in [0,T]$ and $x\in \mathbb{R}^{d}$, the space ${\cal{HP}}$ can
be described as the space of measurable functions $\varphi (s,x)$,
$s\in [0,T], x\in \mathbb{R}^{d}$ such that $K^{\ast} \varphi  \in
L^{2}([0,T])\otimes {\cal{P}}(\mathbb{R}^{d})$, where
$$K^{\ast} \varphi (s,x) = K(T,s) \varphi (s,x) + \int_{0}^{T} \int_{s}^{T}
\left( (\varphi (r,x)-\varphi (s,x) \right) \partial_{1}K(r,s)dr.$$
So, it is necessary to use the transfer operator $K^{\ast}$ and to
check (in the case of the additive noise) that $K^{\ast}g_{t,x} \in
L^{2}([0,T])\otimes {\cal{P}}(\mathbb{R}^{d})$ where $g_{t,x}(s,y)
=p_{t-s}(x-y)1_[0,t]$ which is in principle rather technical (in the
of the stochastic heat equation with multiplicative
fractional-colored noise, one needs to deal with   the tensor
product operator $(K^{\ast})^{\otimes n}$ which has a complicated
expression).   }
\end{remark}

\begin{remark}
{\rm If $H=1/2$ and $f$ is an arbitrary kernel, it was proved in
\cite{dalang99} (using different methods) that the sufficient
condition for the existence of the solution in $\bR_{+} \times
\bR^d$ of (\ref{heat}) with vanishing initial conditions (i.e.
$u_0=0$), is:
\begin{equation}
\label{nec-suf-cond-adit-noise}
\int_{\bR^d}\frac{\mu(d\xi)}{1+|\xi|^2}<\infty.
\end{equation}
(see Remark 14 of \cite{dalang99}).  When $f$ is the Riesz kernel of
order $\alpha$, or the Bessel kernel of order $\alpha$,
(\ref{nec-suf-cond-adit-noise}) holds if and only if $d<\alpha+2$.
Combining this with Remark \ref{nec-remark}, we conclude that,  in
the case of the two kernels, condition
(\ref{nec-suf-cond-adit-noise}) is also necessary for the existence
of the solution. For an arbitrary kernel $f$, it is not known if
(\ref{nec-suf-cond-adit-noise}) remains a necessary condition for
the existence of the solution.

If $H>1/2$ and $f$ is the Riesz kernel of order $\alpha$, or the
Bessel kernel of order $\alpha$, the necessary and sufficient for
the existence of the solution of the stochastic heat equation with
linear additive noise 
is $d<4H+\alpha$, whereas if $f$ is heat or the Poisson kernel, this
equation has a solution for any $d$ (see \cite{balan-tudor08} and
\cite{balan-tudor08-err}). }
\end{remark}

\section{Relationship with the Local Time}

In this section, we identify a random variable $L_t$, defined
formally as a ``convoluted intersection local time'' of two
independent $d$-dimensional standard Brownian motions, such that
\begin{equation}
\label{n-moment-local-time} \alpha_n(t)=E(L_t^n), \quad \forall n
\geq 1.
\end{equation}

An immediate consequence of (\ref{moment-2-utx}),
(\ref{estimate-alpha-n-prime}) and (\ref{n-moment-local-time}) is
that the second moment of $u_{t,x}$ is bounded by the exponential
moment of $L_t$:
$$E|u_{t,x}|^2 \leq \|u_0\|_{\infty}^2 \sum_{n=0}^{\infty}\frac{1}{n!}
\alpha_n(t) = \|u_0\|_{\infty}^2 \sum_{n=0}^{\infty}\frac{1}{n!}
E(L_t^n)=\|u_0\|_{\infty}^2 E(e^{L_t}),$$ with equality if $u_0=1$.

To show (\ref{n-moment-local-time}), we approximate $\alpha_n(t)$ by
$\{\alpha_{n,\varepsilon}(t)\}_{\varepsilon>0}$, when $\varepsilon
\to 0$, where the constants $\alpha_{n,\varepsilon}(t)$ are chosen
such that:
$$\alpha_{n,\varepsilon}(t)=E(L_{t,\varepsilon}^n), \quad \forall
n \geq 1,$$ for a certain random variable $L_{t,\varepsilon}$.

To identify the approximation constants $\alpha_{n,\varepsilon}(t)$,
we recall the definition (\ref{def-alpha-n}), which says that
$\alpha_n(t)$ is the weighted integral of the function
$\psi^{*(n)}({\bf s}, {\bf t})$. The next lemma gives the exact
calculation for the integrand $\psi^{*(n)}({\bf s}, {\bf t})$.

\begin{lemma}
\label{correction-lemma4-2} We have:
\begin{eqnarray*}
\psi^{*(n)}({\bf s},{\bf t})  & = & (2\pi)^{-nd}\int_{(\bR^d)^n}
\exp\left\{- \frac{1}{2}\sum_{j,k=1}^{n} \sigma_{jk}^* \xi_j \cdot
\xi_k \right\} \mu(d\xi_1) \ldots \mu(d\xi_n),
\end{eqnarray*}
where $\sigma_{jk}^*:=(t-s_j) \wedge (t-s_k)+(t-t_j) \wedge
(t-t_k)$.
\end{lemma}

\begin{remark}
\label{remark-about-correction} {\rm Lemma \ref{correction-lemma4-2}
gives a generalization -and a minor correction- to Lemma 4.2 of
\cite{hu-nualart08}. The correction refers to the fact that the
result of \cite{hu-nualart08} is stated incorrectly with the
constant $\sigma_{jk}=s_j \wedge s_k+t_j \wedge t_k$, instead of
$\sigma_{jk}^*$. However, a trivial change of variables
$s_j':=t-s_j,t_j':=t-t_j$ in the definition (\ref{def-alpha-n}) of
$\alpha_n(t)$ shows that this minor error does not affect the
calculation of $\alpha_n(t)$. We have indeed:
\begin{equation}
\label{definition-alpha} \alpha_n(t)=\left\{
\begin{array}{ll}
\alpha_H^n \int_{[0,t]^{2n}} \prod_{j=1}^{n}|s_j-t_j|^{2H-2}
\psi^{(n)}({\bf s}, {\bf t}) d{\bf s} d {\bf t} & \mbox{if $H>1/2$} \\
\\
\int_{[0,t]^{n}} \psi^{(n)}({\bf s}, {\bf s}) d{\bf s}  & \mbox{if
$H=1/2$}
\end{array} \right.
\end{equation}
 with
\begin{eqnarray*}
\psi^{(n)}({\bf s}, {\bf t})&=& (2\pi)^{-nd}\int_{(\bR^d)^n}
\exp\left\{- \frac{1}{2}\sum_{j,k=1}^{n} \sigma_{jk} \xi_j \cdot
\xi_k \right\} \mu(d\xi_1) \ldots \mu(d\xi_n) \\
&=& \psi^{*(n)}(t{\bf 1}-{\bf s},t{\bf 1}-{\bf t}), \quad
\mbox{where} \ {\bf 1}=(1, \ldots,1) \in \bR^n.
\end{eqnarray*} }
\end{remark}

\noindent {\bf Proof of Lemma \ref{correction-lemma4-2}:} Note that
$$\langle \varphi, \psi \rangle_{\cP(\bR^d)^{\otimes n}}=(2\pi)^{-nd}
\int_{(\bR^d)^n} \cF \varphi(\xi_1,\ldots, \xi_n) \overline{\cF
\psi(\xi_1, \ldots, \xi_n)} \mu(d\xi_1) \ldots \mu(d\xi_n),$$ where
$\cF$ denotes the Fourier transform. Hence,
\begin{eqnarray*}
\psi^{*(n)}({\bf s},{\bf t}) & = &\langle g_{\bf s}^{(n)}, g_{\bf
t}^{(n)} \rangle_{\cP(\bR^d)^{\otimes n}} \\
& = & (2\pi)^{-nd} \int_{(\bR^d)^n} \cF g_{\bf
s}^{(n)}(\xi_1,\ldots, \xi_n) \overline{\cF g_{\bf t}^{(n)}(\xi_1,
\ldots, \xi_n)} \mu(d\xi_1) \ldots \mu(d\xi_n).
\end{eqnarray*}

It was shown in the proof of Lemma 4.2 of \cite{hu-nualart08} that:
\begin{eqnarray*}
\cF g_{\bf s}(\xi_1, \ldots, \xi_n) &=& E\left[\prod_{j=1}^{n}
e^{-i \xi_j \cdot [x-(B_{t}^{1}-B_{s_j}^{1})]} \right] \\
\cF g_{\bf t}(\xi_1, \ldots, \xi_n) &=& E\left[\prod_{j=1}^{n} e^{-i
\xi_j \cdot [x-(B_{t}^{2}-B_{t_j}^{2})]} \right],
\end{eqnarray*}
where $B^1=(B_t^1)_{t \geq 0}$ and $B^2=(B_t^2)_{t \geq 0}$ are
independent $d$-dimensional standard Brownian motions. Hence,
\begin{eqnarray*}
\psi^{*(n)} ({\bf s}, {\bf t}) &=& (2 \pi)^{-nd} \int_{(\bR^d)^n} E
\left[\prod_{j=1}^{n} e^{-i \xi_j \cdot
[(B_{s_j}^1-B_t^1)-(B_{t_j}^2-B_t^2)]} \right]
\mu(d\xi_1) \ldots \mu(d \xi_n) \\
\end{eqnarray*}

We begin to evaluate the integrand of the above integral. We denote
$\xi_{j}=(\xi_{j,1}, \ldots, \xi_{j,d})$, $B_t^1=(B_{t,1}^1, \ldots,
B_{t,d}^1)$ and $B_t^2=(B_{t,1}^2, \ldots, B_{t,d}^2)$. We observe
that for any $j=1, \ldots, n$ fixed, the random variables
$$(B_{s_j,l}^{1}-B_{t,l}^{1})-(B_{t_j,l}^{2}-B_{t,l}^{2}), \quad
l=1, \ldots,d, \quad \mbox{are i.i.d.},$$ with the same distribution
as $(b_{s_j}^{1}-b_{t}^{1})-(b_{t_j}^{2}-b_{t}^{2})$, where
$b^1=(b_t^1)_{t \geq 0}$ and $b^2=(b_t^2)_{t \geq 0}$ are
independent $1$-dimensional standard Brownian motions. Hence,
\begin{eqnarray*}
\lefteqn{E \left[\prod_{j=1}^{n} e^{-i \xi_j \cdot
[(B_{s_j}^{1}-B_{t}^{1})-(B_{t_j}^{2}-B_{t}^{2})]} \right] =
\prod_{l=1}^{d}  E\left[ \prod_{j=1}^{n} e^{-i \xi_{j,l}
[(B_{s_j,l}^{1}-B_{t,l}^{1})-(B_{t_j,l}^{2}-B_{t,l}^{2})]} \right]
} \\
& & = \prod_{l=1}^{d}  E\left[ \prod_{j=1}^{n} e^{-i \xi_{j,l}
[(b_{s_j}^{1}-b_{t}^{1})-(b_{t_j}^{2}-b_{t}^{2})]} \right] =
\prod_{l=1}^{d} \exp \left\{-\frac{1}{2} \sum_{j,k=1}^{n}
\sigma_{jk}^{*} \xi_{j,l} \xi_{k,l}  \right\} \\
& & = \exp \left\{-\frac{1}{2}\sum_{j,k=1}^{n} \sigma_{jk}^* \xi_{j}
\cdot \xi_k \right\},
\end{eqnarray*}
where for the second last equality we used the fact that the vector
$$((b_{s_1}^{1}-b_{t}^{1})-(b_{t_1}^2-b_{t}^{2}), \ldots,
(b_{s_n}^{1}-b_{t}^{1})-(b_{t_n}^{2}-b_{t}^{2}))$$ has a normal
distribution with mean zero and covariance matrix
$(\sigma_{jk}^*)_{1 \leq j,k \leq n}$. This concludes the proof of
the lemma. $\Box$

\vspace{3mm}

In what follows, we use the alternative definition
(\ref{definition-alpha}) of $\alpha_n(t)$, given in Remark
\ref{remark-about-correction}. The idea is to find a suitable
approximation for the integrand $\psi^{(n)}({\bf s}, {\bf t})$, by
replacing the Dirac function $\delta_0(x)$ with the heat kernel
$p_{\varepsilon}(x)$. This approximation turns out to be:
$$\psi_{\varepsilon}^{(n)}({\bf s}, {\bf t}):=E \left[
\int_{(\bR^d)^n} p_{\varepsilon}(B_{s_1}^{1}-B_{t_1}^{2}-y_1) \ldots
p_{\varepsilon} (B_{s_n}^{1}-B_{t_n}^{2}-y_n)f(y_1) \ldots f(y_n) d
{\bf y}\right],$$ where $B^1=(B_t^1)_{t \geq 0}$ and $B^2=(B_t^2)_{t
\geq 0}$ are independent $d$-dimensional standard Brownian motions,
and we denote ${\bf y}=(y_1, \ldots,y_n)$.

More precisely, we have the following result.

\begin{lemma}
\label{lemma-form-psi-epsilon} Suppose that $\mu(d \xi)=g(\xi)d
\xi$, i.e. $f=\cF g$. Then
$$\psi_{\varepsilon}^{(n)}({\bf s}, {\bf t})=(2 \pi)^{-nd}
\int_{(\bR^d)^n} \exp \left\{-\frac{1}{2} \sum_{j,k=1}^{n}
\sigma_{jk} \xi_j \cdot \xi_k -\frac{\varepsilon}{2}\sum_{j=1}^{n}
|\xi_j|^2 \right\} \mu(d \xi_1) \ldots \mu(d \xi_n).$$
\end{lemma}

\noindent {\bf Proof:} We first calculate the inverse Fourier
transform of $p_{\varepsilon}*f$:
$$\cF^{-1}(p_{\varepsilon}*f)(\xi)=\cF^{-1} p_{\varepsilon}(\xi) \
\cF^{-1} f(\xi)=(2 \pi)^{-d}e^{-\varepsilon|\xi|^2/2}g(\xi).$$ This
shows that $(p_{\varepsilon}*f) (x)=(2 \pi)^{-d} \cF [
e^{-\varepsilon|\xi|^2/2}g(\xi)](x)$, i.e.
\begin{equation}
\label{Fourier-p*f} \int_{\bR^d}\frac{1}{(2 \pi \varepsilon)^{d/2}}
e^{-|x-y|^2/(2 \varepsilon)}f(y)dy = (2 \pi)^{-d}\int_{\bR^d} e^{-i
\xi \cdot x} e^{-\varepsilon |\xi|^2/2}g(\xi)d \xi.
\end{equation}

Using (\ref{Fourier-p*f}) with $x=B_{s_j}^1-B_{t_j}^2$, we obtain:
$$\int_{\bR^d}\frac{1}{(2 \pi \varepsilon)^{d/2}}
e^{-|B_{s_j}^1-B_{t_j}^2-y_j|^2/(2 \varepsilon)}f(y_j)dy_j = (2
\pi)^{-d}\int_{\bR^d} e^{-i \xi_j \cdot (B_{s_j}^1-B_{t_j}^2)}
e^{-\varepsilon|\xi_j|^2/2}g(\xi_j)d \xi_j.$$

Therefore,
\begin{eqnarray*}
\psi_{\varepsilon}^{(n)}({\bf s}, {\bf t}) &=& E \left[
\prod_{j=1}^{n} \int_{\bR^d}\frac{1}{(2 \pi \varepsilon)^{d/2}}
e^{-|B_{s_j}^1-B_{t_j}^2-y_j|^2/(2\varepsilon)}f(y_j) dy_j \right]
\\
&=& (2 \pi)^{-nd} E \left[ \prod_{j=1}^{n} \int_{\bR^d} e^{-i \xi_j
\cdot (B_{s_j}^1-B_{t_j}^2)} e^{-\varepsilon|\xi_j|^2/2}g(\xi_j)d
\xi_j
\right] \\
&=& (2 \pi)^{-nd} \int_{(\bR^d)^n} E \left[\prod_{j=1}^{n}  e^{-i
\xi_j \cdot (B_{s_j}^1-B_{t_j}^2)} \right] e^{-\varepsilon
\sum_{j=1}^{n}|\xi_j|^2/2} \mu(d \xi_1) \ldots \mu(d \xi_n) \\
&=& (2 \pi)^{-nd} \int_{(\bR^d)^n} e^{-\sum_{j,k=1}^{n} \sigma_{jk}
\xi_j \cdot \xi_k/2} e^{-\varepsilon \sum_{j=1}^{n}|\xi_j|^2/2}
\mu(d \xi_1) \ldots \mu(d \xi_n).
\end{eqnarray*}
$\Box$

\vspace{3mm}

\begin{remark}
\label{cont-of-g} {\rm Note that the function $h=p_{\varepsilon}*f$
is continuous. To see this, let $(x_n)_n \subset \bR^d$ such that
$x_n \to x$. By (\ref{Fourier-p*f}) and the dominated convergence
theorem, it follows that $h(x_n) \to h(x)$. To justify this, note
that $|e^{-i \xi \cdot x_n} e^{-\varepsilon |\xi|^2/2}g(\xi)| \leq
e^{-\varepsilon |\xi|^2/2}g(\xi)$ for any $\xi \in \bR^d, n \geq 1$
and
$$(2\pi)^{-2d}\int_{\bR^d} e^{-\varepsilon |\xi|^2/2}g(\xi)d\xi=
\int_{\bR^d}|\cF p_{\varepsilon/2}(\xi)|^2
g(\xi)d\xi=\|p_{\varepsilon/2}\|_{\cP(\bR^d)}^2<\infty,$$ since
$p_{\varepsilon/2} \in L^2(\bR^d) \subset \cP(\bR^d)$. }
\end{remark}

We are now ready to define the approximation constants
$\alpha_{n,\varepsilon}(t)$:
$$\alpha_{n,\varepsilon}(t)= \left\{
\begin{array}{ll} \alpha_H^n \int_{[0,t]^{2n}}\prod_{j=1}^{n} |t_j-s_j|^{2H-2}
\psi_{\varepsilon}^{(n)}({\bf s}, {\bf t}) d{\bf s} d{\bf t}
& \mbox{if $H>1/2$} \\
\\
\int_{[0,t]^{n}} \psi_{\varepsilon}^{(n)}({\bf s}, {\bf s}) d{\bf s}
& \mbox{if $H=1/2$}
\end{array} \right. .
$$
Note that $\alpha_{n,\varepsilon}(t)=E(L_{t,\varepsilon}^{n})$,
where
$$L_{t,\varepsilon}:=\left\{
\begin{array}{ll}
\alpha_H \int_{0}^{t} \int_{0}^{t} \int_{\bR^d}
|r-s|^{2H-2}p_{\varepsilon}(B_{r}^{1}-B_{s}^{2}-y)f(y)dy dr ds & \mbox{if $H>1/2$} \\
\\
\\
\int_{0}^{t} \int_{\bR^d}
p_{\varepsilon}(B_{s}^{1}-B_{s}^{2}-y)f(y)dy ds  & \mbox{if $H=1/2$}
\end{array} \right.$$

\noindent The random variable $L_{t,\varepsilon}$ is an
approximation of the ``convoluted intersection local time'' $L_{t}$,
written formally as:
$$L_t =\left\{
\begin{array}{ll}
\alpha_H  \int_{0}^{t} \int_{0}^{t} \int_{\bR^d}
|r-s|^{2H-2}\delta_{0}(B_{r}^{1}-B_{s}^{2}-y)f(y)dy dr ds & \mbox{if $H>1/2$} \\
\\
\int_{0}^{t} \int_{\bR^d} \delta_{0}(B_{s}^{1}-B_{s}^{2}-y)f(y)dy ds
& \mbox{if $H=1/2$}
\end{array} \right.$$

\begin{remark}
{\rm We mention that this approximation procedure has been
intensively used in several papers dealing with the chaos expansion
of the local time and Tanaka's formulas for Brownian motion
 (see e.g. \cite{NV}) or fractional Brownian motion (see \cite{CNT}).
 Recall that  the local time of the Brownian motion can be formally
 written as $L(t,x)=\int_{0}^{t} \delta_{0} (B_{s}-x)ds$ where $\delta_{0} $
 is the delta Dirac function. Usually, to obtain the chaos  expansion of
  $L(t,x)$ one approximates $\delta _{0}(B_{s}-x)$ by the Gaussian
  kernel $p_{\epsilon}$.}
\end{remark}

More generally (and for the sake of a result encountered later in
the sequel), if $\eta:[0,t]^2 \to \bR_{+}$ is an arbitrary function
such that $\eta(r,s)=\eta(t-r,t-s)$ for all $r,s \in [0,t]$, we
define
$$L_{t,\varepsilon}(\eta):=\int_0^t \int_0^t \int_{\bR^d}\eta(r,s)
p_{\varepsilon}(B_{r}^{1}-B_{s}^{2}-y)f(y)dy dr ds.$$ Then
$\alpha_{n,\varepsilon}(t,\eta)=E(L_{t,\varepsilon}(\eta)^n)$, where
$\alpha_{n,\varepsilon}(t,\eta):= \int_{[0,t]^{2n}}\prod_{j=1}^{n}
\eta(s_j,t_j) \psi_{\varepsilon}^{(n)}({\bf s}, {\bf t}) d{\bf s}
d{\bf t}$. Let
\begin{equation}
\label{def-alpha-n-eta} \alpha_{n}(t,\eta):=
\int_{[0,t]^{2n}}\prod_{j=1}^{n} \eta(s_j,t_j) \psi^{(n)}({\bf s},
{\bf t}) d{\bf s} d{\bf t},
\end{equation}
and note that $\alpha_{n, \varepsilon}(t,\eta) \leq
\alpha_{n}(t,\eta)$ for all $\varepsilon>0$. Note that $\alpha_{n,
\varepsilon}(t,\eta) \leq \alpha_{n, \varepsilon'}(t,\eta)$ if
$0<\varepsilon'<\varepsilon$.

\begin{lemma}
\label{existence-of-L} Let $t>0$ be arbitrary. a) If
$\alpha_2(t,\eta)<\infty$, then
\begin{equation}
\label{convergence-L_delta-epsilon} \lim_{\varepsilon,\delta
\downarrow
0}E(L_{t,\varepsilon}(\eta)L_{t,\delta}(\eta))=\alpha_2(t,\eta),
\end{equation} and there exists a random variable
$L_{t}(\eta):=\lim_{\varepsilon \downarrow
0}L_{t,\varepsilon}(\eta)$ in $L^{2}(\Omega)$.

b) If $\alpha_n(t,\eta)<\infty$ for all $n \geq 1$, then the random
variable $L_t$, defined in part a), is $p$-integrable for any $p
\geq 2$, and
\begin{equation}
\label{Lp-convergence-L_delta-epsilon} \lim_{\varepsilon \downarrow
0}E|L_{t,\varepsilon}(\eta)-L_{t}(\eta)|^p=0, \quad \mbox{for all} \
p \geq 2.
\end{equation}
In particular, $E(L_t(\eta)^n)=\lim_{\varepsilon \downarrow
0}E(L_{t,\varepsilon}(\eta)^n)=\alpha_n(t,\eta)$ for all $n \geq 1$.
\end{lemma}

The random variable $L_t(\eta)$ defined in Lemma
\ref{existence-of-L} depends on $B^1$ and $B^2$, and could be
denoted by:
$$L_{t}^{B^1,B^2}(\eta)=\int_{0}^{t} \int_{0}^{t} \int_{\bR^d}
\eta(r,s)\delta_{0}(B_{r}^{1}-B_{s}^{2}-y)f(y)dy dr ds.$$ This
notation emphasizes dependence on $B^1,B^2$, and the formal
interpretation of $L_t(\eta)$ as a ``convoluted intersection local
time'' of $B^1$ and $B^2$.

\vspace{3mm}

\noindent {\bf Proof of Lemma \ref{existence-of-L}:} As in
\cite{hu-nualart08}, the proof follows by classical methods. We
include it for the sake of completeness. To simplify the writing, we
omit $\eta$ in the arguments below.

a) Note that $E(L_{t,\varepsilon}L_{t,\delta})=\int_{[0,t]^{4}}
\eta(s_1,t_1) \eta(s_2,t_2) \psi_{\varepsilon,\delta}^{(2)}({\bf s},
{\bf t}) d{\bf s} d{\bf t}$, where
\begin{eqnarray*}
\psi_{\varepsilon,\delta}^{(2)}({\bf s}, {\bf t}) & := & E
\left[\int_{(\bR^d)^2} p_{\varepsilon}(B_{s_1}^{1}-B_{t_1}^{2}-y_1)
p_{\varepsilon}(B_{s_2}^{1}-B_{t_2}^{2}-y_2) f(y_1)f(y_2)dy_1 dy_2
\right] \\
& = & (2\pi)^{-2d} \int_{(\bR^d)^2} \exp\left\{-\frac{1}{2}
\sum_{j,k=1}^{2} \sigma_{jk}\xi_j \cdot \xi_k
-\frac{\varepsilon}{2}|\xi_1|^2 -\frac{\delta}{2}|\xi_2|^2 \right\}
\mu(d\xi_1) \mu(d\xi_2).
\end{eqnarray*}
(The second equality above can be proved using the same argument as
in the proof of Lemma \ref{lemma-form-psi-epsilon}.) Then
$\lim_{\varepsilon, \delta \downarrow
0}\psi_{\varepsilon,\delta}^{(2)} ({\bf s}, {\bf t})=\psi^{(2)}({\bf
s}, {\bf t})$. Relation (\ref{convergence-L_delta-epsilon}) follows
by the dominated convergence theorem, since
$\psi_{\varepsilon,\delta}^{(2)} ({\bf s}, {\bf t}) \leq \psi^
{(2)}({\bf s}, {\bf t})$ for all $\varepsilon, \delta
>0$, and
$$\int_{[0,t]^{4}} \eta(s_1,t_1) \eta(s_2,t_2) \psi^{(2)}({\bf s},
{\bf t})d{\bf s} d{\bf t}=\alpha_2(t)<\infty.$$

>From here, we also infer that $\lim_{\varepsilon \downarrow
0}E(L_{t,\varepsilon}^2)=\alpha_2(t)$, and hence
\begin{equation}
\label{conv-e-d} \lim_{\varepsilon,\delta \downarrow
0}E|L_{t,\varepsilon}-L_{t,\delta}|^2=\lim_{\varepsilon \downarrow
0}E(L_{t,\varepsilon}^2)+\lim_{\delta \downarrow
0}E(L_{t,\delta}^2)-2\lim_{\varepsilon,\delta \downarrow
0}E(L_{t,\varepsilon}L_{t,\delta})=0.
\end{equation}

Let $(\varepsilon_n)_n \downarrow 0$ be arbitrary. From
(\ref{conv-e-d}), it follows that $(L_{t,\varepsilon_n})_{n}$ is a
Cauchy sequence in $L^{2}(\Omega)$. Hence, there exists $L_t \in
L^2(\Omega)$ such that $E|L_{t,\varepsilon_n}-L_t|^2 \to 0$. If
$(\varepsilon_n')_n \downarrow 0$ is another sequence and
$E|L_{t,\varepsilon_n'}-L_t'|^2 \to 0$ for some $L_t' \in
L^2(\Omega)$, then $E|L_t-L_t'|^2 \leq
E|L_t-L_{t,\varepsilon_n}|^2+E|L_{t,\varepsilon_n}-L_{t,\varepsilon_n'}|^2+
E|L_{t,\varepsilon_n'}-L_t'|^2 \to 0$, i.e. $E|L_t-L_{t}'|^2=0$.
This shows that $L_t$ does not depend on $(\varepsilon_n)_n$.

b) Let $p \geq 2$ be fixed. Let $(\varepsilon_n)_n \downarrow 0$ be
arbitrary. We will prove that $E|L_{t,\varepsilon_n}-L_t|^p \to 0$,
by using the fact that, in a metric space, $x_n \to x$ if and only
if for any subsequence $N' \subset \bN$ there exists a
sub-subsequence $N'' \subset N'$ such that $x_{n} \to x$, as $n \to
\infty, n \in N''$ (see e.g. p.15 of \cite{billingsley68}).

Let $N' \subset \bN$ be an arbitrary subsequence. By part a), as $n
\to \infty,n \in N'$, $L_{t,\varepsilon_{n}} \to L_t$ in
$L^{2}(\Omega)$. Hence, $L_{t,\varepsilon_{n}} \to L_t$ in
probability, and there exists a sub-subsequence $N'' \subset N'$
such that $L_{t,\varepsilon_{n}} \to L_t$ a.s., as $n \to \infty, n
\in N''$. Note that $(L_{t,\varepsilon})_{\varepsilon>0}$ is
uniformly integrable, since
$$\sup_{\varepsilon>0}E(L_{t,\varepsilon}^n)=
\sup_{\varepsilon>0}\alpha_{n,\varepsilon}(t) \leq
\alpha_n(t)<\infty, \quad \mbox{for} \ n \geq 2.$$

\noindent By Theorem 16.14 of \cite{billingsley95}, it follows that
$|L_t|^p$ is integrable and $E|L_{t,\varepsilon_{n}}-L_t|^p \to 0$,
as $n \to \infty, n \in N''$. $\Box$

\vspace{3mm}

The next two results are the analogues of Propositions 3.1 and 3.2
of \cite{hu-nualart08} in the case of a colored noise. We denote
$\Phi(x,a)=\sum_{n=0}^{\infty}x^n/(n!)^a$ for $x >0$ and $a \geq 0$.
Note that $\Phi(x,a)<\infty$ if and only if $a>0,x>0$ or $a=0, x \in
(0,1)$.

\begin{proposition}
\label{exp-L-bound-2} Suppose that $\eta :[0,t]^2 \to \bR_+$
satisfies the following condition:
\begin{equation}
\label{cond-C1} \|\eta\|_{1,t}:=\max \left( \sup_{s \in [0,t]}
\int_0^t \eta(r,s)dr, \sup_{r \in [0,t]}\int_0^t \eta(r,s)ds
\right)<\infty.
\end{equation}

(i) If $f$ is the Riesz kernel of order $\alpha$, or the Bessel
kernel of order $\alpha<d$, and $d<2+\alpha$, then
$\lim_{\varepsilon \downarrow 0}L_{t,\varepsilon}(\eta)=L_{t}(\eta)$
exists in $L^p(\Omega)$ for all $p \geq 2$, and $$\sup_{\varepsilon
>0}E\left[\exp\left(\lambda L_{t,\varepsilon}(\eta)\right) \right]
\leq C_{\alpha,d}^{*} \Phi\left(\lambda
D(t),1-\frac{d-\alpha}{2}\right), \quad \mbox{for all} \
\lambda>0,$$ where $C_{\alpha,d}^*$ is a constant depending on
$\alpha$ and $d$, and
$$D(t)=D_{\alpha,d} 2^{-(d-\alpha)/2} \|\eta\|_{1,t} \Gamma
\left(1-\frac{d-\alpha}{2}\right)
\left(1-\frac{d-\alpha}{2}\right)^{-[1-(d-\alpha)/2]}
t^{1-(d-\alpha)/2}.$$

(ii) If  $f$ is the heat kernel of order $\alpha$, or the Poisson
kernel of order $\alpha$, then $\lim_{\varepsilon \downarrow
0}L_{t,\varepsilon}(\eta)=L_{t}(\eta)$ exists in $L^p(\Omega)$, for
all $p \geq 2$, and $$\sup_{\varepsilon
>0}E\left[\exp\left(\lambda L_{t,\varepsilon}(\eta)\right) \right]
\leq \exp (\lambda C(t)), \quad \mbox{for all} \ \lambda>0,$$ where
$C(t)=C_{\alpha,d} \|\eta\|_{1,t} t$.
\end{proposition}

\noindent {\bf Proof:} We use the definition (\ref{def-alpha-n-eta})
of $\alpha_n(t,\eta)$, and Lemma \ref{estimate-psi} for the
estimation of $\psi^{*(n)}({\bf s},{\bf t})$.

(i) Let $h=-(d-\alpha)/2$. Using the Cauchy-Schwartz inequality,
condition (\ref{cond-C1}), and Lemma \ref{calcul-In}, we get:
\begin{eqnarray*}
\alpha_n(t,\eta) & \leq & \left(D_{\alpha,d} 2^{-(d-\alpha)/2}
\right)^n \int_{[0,t]^{2n}} \prod_{j=1}^{n}\eta(s_j,t_j)[\beta({\bf
s})
\beta({\bf t})]^{-(d-\alpha)/4}d{\bf s}d{\bf t} \\
 & \leq &
\left(D_{\alpha,d} 2^{-(d-\alpha)/2} \right)^n \int_{[0,t]^{2n}}
\prod_{j=1}^{n}\eta(s_j,t_j) [\beta({\bf
s})]^{-(d-\alpha)/2} d{\bf s}d{\bf t} \\
& \leq & \left(D_{\alpha,d} 2^{-(d-\alpha)/2} \right)^n
\|\eta\|_{1,t}^n
\int_{[0,t]^{n}} [\beta({\bf s})]^{-(d-\alpha)/2} d{\bf s} \\
&=& \left(D_{\alpha,d} 2^{-(d-\alpha)/2} \|\eta\|_{1,t} \right)^n n!
\ I_n(t,h) \\
& = & \Gamma(1+h)\left(D_{\alpha,d} 2^{-(d-\alpha)/2} \|\eta\|_{1,t}
\Gamma(1+h) t^{1+h}\right)^n \frac{n!}{\Gamma(n(1+h)+1)} \\
& \leq & C_{\alpha,d}^* D(t)^n (n!)^{-h}.
\end{eqnarray*} (The last
inequality follows by relation (3.19) of \cite{hu-nualart08}.) The
first statement follows by Lemma \ref{existence-of-L}. The second
statement follows since,
$$
E[e^{\lambda L_{t,\varepsilon}(\eta)}] =\sum_{n=0}^{\infty}
\frac{\lambda^n}{n!} \alpha_{n,\varepsilon}(t,\eta) \leq
\sum_{n=0}^{\infty} \frac{\lambda^n}{n!} \alpha_{n}(t,\eta)  \leq
C_{\alpha,d}^* \sum_{n=0}^{\infty}\frac{[\lambda
D(t)]^n}{(n!)^{1+h}},$$ and the last sum is finite for all
$\lambda>0$, since $1+h>0$.

(ii) In this case,
\begin{eqnarray*}
\alpha_n(t,\eta) & \leq C_{\alpha,d}^{n}\int_{[0,t]^{2n}}
\prod_{j=1}^{n}\eta(s_j,t_j)d{\bf s}d{\bf t} \leq C_{\alpha,d}^n
\|\eta\|_{1,t}^n t^n=C(t)^n,
\end{eqnarray*}
and
$$E[e^{\lambda L_{t,\varepsilon}(\eta)}]
\leq \sum_{n=0}^{\infty} \frac{\lambda^n}{n!} \alpha_{n}(t,\eta)
\leq \sum_{n=0}^{\infty}\frac{[\lambda C(t)]^n}{n!}=e^{\lambda
C(t)}.$$
 $\Box$

\begin{proposition}
\label{exp-L-bound} Suppose that $\eta :[0,t]^2 \to \bR_+$ satisfies
the following condition: there exist $\gamma>0$ and $1/2<H<1$, such
that
\begin{equation}
\label{cond-C2} \eta(r,s) \leq \gamma |r-s|^{2H-2}, \quad \forall
r,s \in [0,t].
\end{equation}

(i) If $f$ is the Riesz kernel of order $\alpha$, or $f$ is the
Bessel kernel of order $\alpha<d$, and $d \leq 2+\alpha$, then
$\lim_{\varepsilon \downarrow 0}L_{t,\varepsilon}(\eta)=L_{t}(\eta)$
exists in $L^p(\Omega)$, for all $p \geq 2$, and
$$
\sup_{\varepsilon>0}E\left[\exp\left(\lambda
L_{t,\varepsilon}(\eta)\right)
\right]\leq C_{H,d,\alpha}^* \Phi\left(\lambda D(t),1-\frac{d-\alpha}{2}\right), 
\quad \mbox{for all} \ 0<\lambda< \lambda_0(t),$$
 where $C_{H,d,\alpha}^*$ is a constant depending on $H,d$ and $\alpha$,
$$D(t)= D_{\alpha,d} 2^{-(d-\alpha)/2} \frac{\gamma}{\alpha_H}
\beta_H^{2}\Gamma \left(1-
 \frac{d-\alpha}{4H}\right)^{2H} \left(1-\frac{d-\alpha}{4H}
 \right)^{-[2H-(d-\alpha)/(2)]}t^{2H-(d-\alpha)/2}$$ and
$$\lambda_0(t)=\left\{
\begin{array}{ll}  \left(1-\frac{1}{2H} \right)^{2H-1}D_{\alpha,d}^{-1}2
\gamma^{-1} \beta_H^{-2} \Gamma \left(1- \frac{1}{2H}
\right)^{-2H}t^{1-2H}
& \mbox{if $d=2+\alpha$} \\
\infty & \mbox{if $d<2+\alpha$}
\end{array} \right.$$

(ii) If $f$ is the heat kernel of order $\alpha$ or the Poisson
kernel of order $\alpha$, then then $\lim_{\varepsilon \downarrow
0}L_{t,\varepsilon}(\eta)=L_{t}(\eta)$ exists in $L^p(\Omega)$, for
all $p \geq 2$, and $$\sup_{\varepsilon>0}E\left[\exp\left(\lambda
L_{t,\varepsilon}(\eta)\right) \right]\leq \exp (\lambda C(t)),
\quad \mbox{for all} \ \lambda>0,$$
 where $C(t)=C_{\alpha,d}\gamma
t^{2H}/\alpha_H$.
\end{proposition}

\noindent {\bf Proof:} The proof is similarly to Proposition
\ref{estimate-alpha-n}. We use the definition
(\ref{def-alpha-n-eta}) of $\alpha_n(t,\eta)$, Lemma
\ref{estimate-psi} and condition (\ref{cond-C2}).

(i) We have: $$\alpha_n(t,\eta) \leq C_{H,d,\alpha}^* D(t)^{n} \
(n!)^{(d-\alpha)/2}<\infty.$$

\noindent The first statement follows by Lemma \ref{existence-of-L}.
The other statement follows, since
$$E[e^{\lambda L_{t,\varepsilon}(\eta)}] =
\sum_{n=0}^{\infty}\frac{\lambda^n}{n!}\alpha_{n,\varepsilon}(t,\eta)
\leq \sum_{n=0}^{\infty}\frac{\lambda^n}{n!}\alpha_{n}(t,\eta)\leq
C_{H,d,\alpha}^* \sum_{n=0}^{\infty}\frac{[\lambda
D(t)]^n}{(n!)^{1-(d-\alpha)/2}}.$$

\noindent If $d-\alpha=2$, then the last sum is finite for all
$0<\lambda<\lambda_0(t):=1/D(t)$. If $d-\alpha<2$, then then last
sum is finite for all $\lambda>0$.

(ii) The result follows, since:
\begin{eqnarray*}
\alpha_{n}(t,\eta) & \leq & C_{\alpha,d}^n \int_{[0,t]^{2n}}
\prod_{j=1}^{n}\eta(s_j,t_j)d{\bf s} d{\bf t} \leq C_{\alpha,d}^n
\gamma^n \int_{[0,t]^{2n}} \prod_{j=1}^{n}|s_j-t_j|^{2H-2}d{\bf s}
d{\bf t} \\ &=&  \left( C_{\alpha,d}\frac{\gamma}{\alpha_H}
t^{2H}\right)^n =C(t)^n
\end{eqnarray*}
and hence
$$E[e^{\lambda L_{t,\varepsilon}(\eta)}]
\leq \sum_{n=0}^{\infty}\frac{\lambda^n}{n!}\alpha_{n}(t,\eta)\leq
\sum_{n=0}^{\infty} \frac{[\lambda C(t)]^n}{n!}.$$ $\Box$

We now introduce the approximation technique of \cite{hu-nualart08},
which will yield simultaneously the existence of the solution of
(\ref{heat}) and some representation formulas for the moments of
this solution. We review briefly this powerful technique, which has
been introduced only recently in the literature. The idea is to
smooth the noise $\dot W$, solve the equation driven by the smoothen
noise, and then show that the solution of the ``smoothen'' equation
converges to the solution of (\ref{heat}).

For any $\varepsilon,\delta>0$, let $\varphi_{\delta}(t)=\delta^{-1}
1_{[0,\delta]}(t)$ and
$$\dot W_{t,x}^{\varepsilon,\delta}=\int_0^t \int_{\bR^d}
\varphi_{\delta}(t-s)p_{\varepsilon}(x-y)d W_{s,y}.$$

\noindent Note that the noise $\dot W^{\varepsilon,\delta}$ can be
viewed as a ``mollification'' of $\dot W$, with rate $\delta$ in the
time variable and rate $\sqrt{\varepsilon}$ is the space variable,
since
$$\varphi_{\delta}=\frac{1}{\delta}\varphi\left(\frac{t}{\delta}\right)
\quad \mbox{and} \quad
p_{\varepsilon}(x)=\frac{1}{(\sqrt{\varepsilon})^d}\phi\left(
\frac{x}{\sqrt{\varepsilon}} \right),$$ with
$\varphi(t)=1_{[0,1]}(t)$ and $\phi(x)=(2\pi)^{-d/2}e^{-|x|^2/2}$.
(Recall that the function $u^{(\varepsilon)}$, defined by
$u^{(\varepsilon)} (x)=\int_{\bR^n}\psi_{\varepsilon}(x-y)u(y)dy$,
is a ``mollification'' of the function $u$ on $\bR^n$, if
$\psi_{\varepsilon}(x)=\varepsilon^{-n}\psi(x/\varepsilon)$ and
$\psi \geq 0$ is such that $\int_{\bR^n} \psi(x)dx=1$.) Therefore,
this approximation procedure can be regarded as a stochastic version
of the ``approximation to the identity'' technique, encountered in
the PDE literature.

\vspace{3mm}

We consider the following ``approximation'' of equation
(\ref{heat}):
\begin{eqnarray}
\label{heat-epsilon-delta} \frac{\partial
u^{\varepsilon,\delta}}{\partial t} &=& \frac{1}{2} \Delta
u^{\varepsilon,\delta}
+ u^{\varepsilon,\delta} \dot {W}^{\varepsilon,\delta} , \quad t>0, x\in \bR^{d} \\
\nonumber u_{0,x}^{\varepsilon,\delta} &=& u_0(x), \quad x \in
\bR^d.
\end{eqnarray}

We introduce now the rigorous meaning for the solution of
(\ref{heat-epsilon-delta}), which could be derived formally from the
mild or evolution version of the equation, by applying the
stochastic Fubini theorem.

\begin{definition}
\label{def-sol-ed} An $(\cF_t)_t$-adapted square-integrable process
$u=\{u_{t,x}^{\varepsilon,\delta}; (t,x) \in \bR_{+} \times
\bR^{d}\}$ is a {\bf solution to} (\ref{heat-epsilon-delta}) if for
any $(t,x) \in \bR_{+} \times \bR^d$, the process
$$\left\{Y_{r,z}^{t,x,\varepsilon,\delta}=1_{[0,t]}(r) \int_0^t
\int_{\bR^d}p_{t-s}(x-y)\varphi_{\delta}(s-r)p_{\varepsilon}(y-z)
 u_{s,y}^{\varepsilon,\delta}dy ds; \ (r,z) \in \bR_{+} \times
\bR^d\right\}$$ exists, is integrable with respect to $W$, and
satisfies
$$u_{t,x}^{\varepsilon,\delta}=p_t
u_0(x)+\int_{0}^{\infty}\int_{\bR^d}
Y_{r,z}^{t,x,\varepsilon,\delta}\delta W_{r,z}.$$
\end{definition}

\noindent By (\ref{duality}), the above definition is equivalent to
saying that for any $(t,x) \in \bR_{+} \times \bR^{d}$, the process
$Y^{t,x,\varepsilon,\delta}$ exists, $u_{t,x}^{\varepsilon,\delta}
\in L^{2}(\Omega)$, $u_{t,x}^{\varepsilon,\delta}$ is
$\cF_t$-measurable and
\begin{equation}
\label{def-sol-u-e-d} E(u_{t,x}^{\varepsilon,\delta}F)=E(F)p_t
u_0(x)+ E \langle Y^{t,x,\varepsilon,\delta}, DF \rangle_{\cH \cP},
\quad \forall F \in \bD^{1,2}.
\end{equation}

Before constructing the solution of (\ref{heat-epsilon-delta}), we
mention few words about the notation. If $X$ and $Y$ are random
variables defined on $(\Omega, \cF, P)$, with values in arbitrary
measurable spaces $\cX$, respectively $\cY$, and $h:\cX \times \cY
\to \bR$ is a measurable function, we define the random variable:
$E^{X}[h(X,Y)](\omega)=\int_{\cX} h(x,Y(\omega)) (P \circ
X^{-1})(dx)$. If $X$ and $Y$ are independent, then
\begin{equation}
\label{relation-A} E[E^{X}[h(X,Y)]]=E[h(X,Y)]=E[E[h(X,Y)|X]],
\end{equation}
where $E [\ \cdot \ ]$ denotes the expectation with respect to $P$,
and $E[\ \cdot \ |X]$ denotes the conditional expectation given $X$.
(This result will be used below with $X=B$ and $Y=W$.)

We have the following result.

\begin{proposition}
\label{exist-sol-ed} The process
$u^{\varepsilon,\delta}=\{u_{t,x}^{\varepsilon,\delta}; (t,x) \in
\bR_{+} \times \bR^d\}$ defined by:
\begin{equation}
\label{def-u-epsilon-delta}
u_{t,x}^{\varepsilon,\delta}:=E^{B}\left[u_0(x+B_t)\exp
\left(\int_0^t \int_{\bR^d}A_{r,y}^{\varepsilon,\delta,B}dW_{r,y}
-\frac{1}{2}\|A^{\varepsilon,\delta,B}\|_{\cH \cP}^2 \right)\right]
\end{equation}
is a solution of (\ref{heat-epsilon-delta}), where
$A_{r,y}^{\varepsilon,\delta,B}=\int_{0}^t
\varphi_{\delta}(t-s-r)p_{\varepsilon}(x+B_s-y)ds$, and $B=(B_t)_{t
\geq 0}$ is a $d$-dimensional standard Brownian motion, independent
of $W$.
\end{proposition}

\noindent {\bf Proof:} The argument is similar to the one used in
the proof of Proposition 5.2 of \cite{hu-nualart08}. We include it
in response to the referee's suggestion. To simplify the notation,
we omit writing $\cH \cP$ in $\|\cdot \|_{\cH \cP}$ and $\langle
\cdot, \cdot \rangle_{\cH \cP}$. We also omit writing
$\varepsilon,\delta$ in $A^{\varepsilon,\delta,B}$, i.e. we denote
$A^{\varepsilon,\delta,B}$ by $A^B$.

For every $\varphi \in \cH \cP$, define
$F_{\varphi}=e^{W(\varphi)-\|\varphi\|^2/2}$. Note that
\begin{equation}
\label{expo-moment-1} E(e^{W(\varphi)})=e^{\|\varphi\|^2/2}, \quad
\forall \varphi \in \cH \cP.
\end{equation}

\noindent Since $\{F_{\varphi}; \varphi \in \cH \cP\}$ is dense in
$\bD^{1,2}$ (see e.g. Lemma 1.1.2 of \cite{nualart06}), it suffices
to prove (\ref{def-sol-u-e-d}) for $F=F_{\varphi}$. Define
$S_{t,x}(\varphi)=E(u_{t,x}^{\varepsilon,\delta}F_{\varphi})$. Using
(\ref{relation-A}),
\begin{eqnarray*}
S_{t,x}(\varphi) &=& E[E^{B}[u_{0}(x+B_t) e^{W(A^B) -\|A^B\|^2/2}]
e^{W(\varphi)-\|\varphi\|^2/2}] \\
&=& E[E^{B}[u_{0}(x+B_t) e^{W(A^B+\varphi) -\|A^B+\varphi\|^2/2}
e^{\langle A^B,\varphi \rangle}]] \\
&=& E[E[u_{0}(x+B_t) e^{W(A^B+\varphi) -\|A^B+\varphi\|^2/2}
e^{\langle A^B,\varphi \rangle}|B]].
\end{eqnarray*}

\noindent Let
$$h(B,W)=u_{0}(x+B_t) e^{W(A^B+\varphi) -\|A^B+\varphi\|^2/2} e^{\langle
A^B,\varphi \rangle}.$$

\noindent Since $B$ and $W$ are independent, $E[h(B,W)|B]=f(B)$,
where
\begin{eqnarray*}
f(b)&=& E[h(b,W)]= E[u_{0}(x+b_t) e^{W(A^b+\varphi)
-\|A^b+\varphi\|^2/2} e^{\langle A^b,\varphi \rangle}]\\
&=& u_{0}(x+b_t) e^{\langle A^b,\varphi \rangle} E[e^{W(A^b+\varphi)
-\|A^b+\varphi\|^2/2} ] \\
&=& u_{0}(x+b_t) e^{\langle A^b,\varphi \rangle}, \quad \mbox{for
any} \ b=(b_t)_{t \geq 0} \in C([0,\infty),\bR^d),
\end{eqnarray*}
where $C([0,\infty),\bR^d)$ denotes the space of continuous
functions $x:[0,\infty) \to \bR^d$, and we used
(\ref{expo-moment-1}) for the last equality. Hence
$E[h(B,W)|B]=u_{0}(x+B_t) e^{\langle A^B,\varphi \rangle}$ and
$$S_{t,x}(\varphi)=E[E[h(B,W)|B]]=E[u_{0}(x+B_t) e^{\langle A^B,\varphi
\rangle}].$$

\noindent By the definition of $A^B$ and Fubini's theorem, we
obtain:
\begin{eqnarray*}
\langle A^B, \varphi \rangle&=& \alpha_H \int_{(\bR_{+} \times
\bR^d)^2} A^B_{r,y}\varphi_{r',y'}|r-r'|^{2H-2} f(y-y')dy dy' dr dr'
\\ &=& \int_{0}^{t}V^{\varepsilon,\delta}(t-s,x+B_s)ds,
\end{eqnarray*}
where $V^{\varepsilon,\delta}(t,x) =\langle
\varphi_{\delta}(t-\cdot) p_{\varepsilon}(x-\cdot),\varphi \rangle$.
Hence:
$$S_{t,x}(\varphi)=E\left[u_0(x+B_t)
\exp \left(\int_0^t V^{\varepsilon,\delta}(t-s,x+B_s)ds \right)
\right].$$

By the Feynman-Kac's formula (see e.g. Theorem 5.7.6 of
\cite{karatzas-shreve91}), $(S_{t,x}(\varphi))_{t,x}$ is a solution
of the Cauchy problem:
\begin{eqnarray*}
\frac{\partial S_{t,x}(\varphi)}{\partial t} &=& \frac{1}{2}\Delta
S_{t,x}(\varphi)+
S_{t,x}(\varphi)V^{\varepsilon,\delta}(t,x), \quad t>0, x \in \bR^d \\
S_{0,x}(\varphi) &=& u_0(x).
\end{eqnarray*}

\noindent Hence,
\begin{eqnarray}
\nonumber S_{t,x}(\varphi)&=& p_y u_0(x)+\int_0^t
\int_{\bR^d}p_{t-s}(x-y)S_{s,y}(\varphi)V^{\varepsilon,\delta}(s,y)dy ds \\
\nonumber &=&  p_y u_0(x)+ \alpha_H E \int_{(\bR_{+} \times
\bR^d)^2} Y_{r,z}^{t,x,\varepsilon,\delta} \varphi(r',z')
F_{\varphi}|r-r'|^{2H-2}f(z-z')dz dz' dr dr' \\
\label{FK-formula} &=& p_y u_0(x)+ E \langle
Y^{t,x,\varepsilon,\delta}, D F_{\varphi} \rangle,
\end{eqnarray}

\noindent where we used Fubini's theorem for the second equality
above and the fact that
$D_{r',z'}F_{\varphi}=\varphi(r',z')F_{\varphi}$ for the third
equality. This concludes the proof of (\ref{def-sol-u-e-d}) for
$F=F_{\varphi}$. $\Box$

\vspace{3mm}

Let $B^i=(B^i_t)_{t \geq 0}, i \geq 1$ be independent
$d$-dimensional standard Brownian motions, independent of $W$.
Suppose that either (\ref{cond-d-alphaf-2}) or
(\ref{cond-d-alphaf-3}) hold. For any pair $(i,j)$ with $i \not =
j$, let $L_{t}^{B^i,B^j}$ be the random variable defined in Lemma
\ref{existence-of-L}, with
$$\eta(r,s)=\left\{
\begin{array}{ll}
\alpha_H|r-s|^{2H-2} & \mbox{if $H>1/2$} \\
1_{\{r=s\}} & \mbox{if $H=1/2$}
\end{array} \right.$$
(By Proposition \ref{estimate-alpha-n}, $\alpha_n(t,\eta)<\infty$
for all $n \geq 1$, and $L_{t}^{B^i,B^j}$ is well-defined.)


\vspace{3mm}

The following result is the main theorem of the present article.

\begin{theorem}
\label{moments} (i) Suppose that $f$ is the Riesz kernel of order
$\alpha$ or the Bessel kernel of order $\alpha<d$, and either
(\ref{cond-d-alphaf-2}) or (\ref{cond-d-alphaf-3}) holds. Then, for
any integer $k \geq 2$, we have:
\begin{equation}
\label{sup-u-epsilon-delta-power-k}
\sup_{\varepsilon,\delta>0}E[(u_{t,x}^{\varepsilon,\delta})^k]<\infty,
\quad \mbox{for all} \ 0<t <t_0(k), x \in \bR^d
\end{equation}
where
$$t_0(k)=\left\{
\begin{array}{ll}
\left[k(k-1)D_{\alpha,d}  2^{-2H}\beta_H^2 \Gamma
\left(1-\frac{1}{2H} \right)^{2H}
\right]^{-1/(2H-1)} & \mbox{if $d=2+\alpha$} \\
\infty & \mbox{if $d<2+\alpha$}
\end{array} \right.$$

For any $0<t<t_0(2)$ and $x \in \bR^d$, the limit
$u_{t,x}:=\lim_{\varepsilon \downarrow 0} \lim_{\delta \downarrow
0}u_{t,x}^{\varepsilon,\delta}$ exists in $L^2(\Omega)$, the process
$u=\{u_{t,x}; (t,x) \in [0,t_0(2)) \times \bR^d\}$ is the unique
solution of (\ref{heat}) in $L^2(\Omega)$, and
$$E[u_{t,x}^2]=E\left[u_{0}(x+B_t^1)u_{0}(x+B_t^2)
\exp \left(L_{t}^{B^1,B^2} \right) \right]:=\gamma_2(t,x).$$

If $x \in \bR^d$ and $t<t_0(M)$ for some $M \geq 3$, then
$\lim_{\varepsilon \downarrow 0} \lim_{\delta \downarrow
0}E|u_{t,x}^{\varepsilon,\delta}-u_{t,x}|^p=0$ for all $2 \leq p
<M$, and for any integer $2 \leq k \leq M-1$,
\begin{equation}
\label{u-power-k} E[u_{t,x}^k]=E\left[\prod_{j=1}^{k}u_{0}(x+B_t^i)
\exp \left(\sum_{1\leq i<j \leq k}L_{t}^{B^i,B^j} \right)
\right]:=\gamma_k(t,x).
\end{equation}

(ii) Suppose that $f$ is the heat kernel of order $\alpha$, or the
Poisson kernel of order $\alpha$. Then the conclusion same as in
part (i) holds, with $t_0(k)=\infty$ for all $k \geq 2$.
\end{theorem}


\noindent {\bf Proof:} The argument is similar to the one used in
the proof of Theorem 5.3 of \cite{hu-nualart08}. At the referee's
request, we include all the details for the reader's convenience. To
ease the exposition, we divide the proof in several steps.

{\em Step 1.} We show that for any integer $k \geq 2$,
\begin{equation}
\label{moment-k-u-epsilon-delta}
E[(u_{t,x}^{\varepsilon,\delta})^k]=E
\left[\prod_{j=1}^{k}u_0(x+B_t^j) \exp \left(\sum_{1 \leq i<j \leq
k} \langle A^{\varepsilon,\delta,B^i}, A^{\varepsilon,\delta,B^j}
\rangle_{\cH \cP} \right) \right].
\end{equation}

By (\ref{def-u-epsilon-delta}), $u_{t,x}^{\varepsilon,\delta}$ can
be expressed as
$$u_{t,x}^{\varepsilon,\delta}=
E^{B^i}\left[u_0(x+B_t^i) \exp \left(\int_0^t \int_{\bR^d}
A_{r,y}^{\varepsilon,\delta,B^i}dW_{r,y}
-\frac{1}{2}\|A^{\varepsilon,\delta,B^i}\|_{\cH \cP}^2 \right)
\right],$$ for any $i=1,\ldots,k$. Taking the product over $i=1,
\ldots,k$ and using the independence of $B^1, \ldots, B^k$, we
obtain that:
\begin{eqnarray*}
(u_{t,x}^{\varepsilon,\delta})^k&=&
\prod_{i=1}^{k}E^{B^i}\left[u_0(x+B_t^i)\exp \left(\int_0^t
\int_{\bR^d}A_{r,y}^{\varepsilon,\delta,B^i}dW_{r,y}
-\frac{1}{2}\|A^{\varepsilon,\delta,B^i}\|_{\cH \cP}^2
\right)\right] \\
&=& E^{B^1, \ldots B^k}\left[\prod_{i=1}^{k} u_0(x+B_t^i)\exp
\left(\int_0^t \int_{\bR^d}A_{r,y}^{\varepsilon,\delta,B^i}dW_{r,y}
-\frac{1}{2}\|A^{\varepsilon,\delta,B^i}\|_{\cH \cP}^2
\right)\right]
\end{eqnarray*}

\noindent Taking the expectation, and using (\ref{relation-A}) with
$X=(B^1, \ldots, B^k):=B$ and $Y=W$, we get:
\begin{eqnarray*}
E[(u_{t,x}^{\varepsilon,\delta})^k]&=& E\left[E
\left[\prod_{i=1}^{k} u_0(x+B_t^i)\exp \left(\int_0^t
\int_{\bR^d}A_{r,y}^{\varepsilon,\delta,B^i}dW_{r,y}
-\frac{1}{2}\|A^{\varepsilon,\delta,B^i}\|_{\cH \cP}^2 \right) | B
\right] \right].
\end{eqnarray*}
Let
$$h(B,W)=\prod_{i=1}^{k} u_0(x+B_t^i)\exp \left(\int_0^t
\int_{\bR^d}A_{r,y}^{\varepsilon,\delta,B^i}dW_{r,y}
-\frac{1}{2}\|A^{\varepsilon,\delta,B^i}\|_{\cH \cP}^2 \right).$$
Then $E[h(B,W)|B]=f(B)$, where
\begin{eqnarray*}
f(b) &=& E[h(b,W)] =E \left[\prod_{i=1}^{k}u_0(x+b_t^i)
e^{W(A^{\varepsilon,\delta,b^i})-\|A^{\varepsilon,\delta,b^i}\|_{\cH
\cP}^2/2} \right] \\
&=& \prod_{i=1}^{k}u_0(x+b_t^i)
e^{-\sum_{i=1}^{k}\|A^{\varepsilon,\delta,b^i}\|_{\cH \cP}^2/2 }
E\left[ e^{W\left(\sum_{i=1}^{k} A^{\varepsilon,\delta,b^i}\right)}\right] \\
&=& \prod_{i=1}^{k}u_0(x+b_t^i)
\exp\left(-\frac{1}{2}\sum_{i=1}^{k}\|A^{\varepsilon,\delta,b^i}\|_{\cH
\cP}^2 +\frac{1}{2}\left\|\sum_{i=1}^{k}
A^{\varepsilon,\delta,b^i}\right\|_{\cH \cP}^2 \right)\\
&=& \prod_{i=1}^{k}u_0(x+b_t^i) \exp \left(\sum_{i<j} \langle
A^{\varepsilon,\delta,b^i}, A^{\varepsilon,\delta,b^j} \rangle_{\cH
\cP} \right)
\end{eqnarray*}
for any $b=(b^1,\ldots,b^k)$ with $b^i=(b_t^i)_{t \geq 0} \in
C([0,\infty),\bR^d)$. (We used (\ref{expo-moment-1}) and the fact
that $W(\varphi+\psi)=W(\varphi)+W(\psi)$ a.s. for any $\varphi,\psi
\in \cH \cP$, which can be checked in $L^2(\Omega)$, using the fact
that $W$ is an isometry between $\cH \cP$ and $L^2(\Omega)$.)
Relation (\ref{moment-k-u-epsilon-delta}) follows, since
$E[(u_{t,x}^{\varepsilon,\delta})^k] = E[E[h(B,W)|B]]=E[f(B)]$.

{\em Step 2.} We prove that for any $(t,x) \in \bR_{+} \times \bR^d$
with $t<t_0(2)$,
\begin{equation}
\label{limit-only-delta}\lim_{\delta \downarrow 0}\langle
A^{\varepsilon,\delta,B^i}, A^{\varepsilon,\delta,B^j} \rangle_{\cH
\cP} = L_{t,2\varepsilon}^{B^i,B^j},    \quad \forall \varepsilon>0,
\quad \forall \omega \in \tilde \Omega_{i,j},
\end{equation}
where $\tilde \Omega_{i,j}=\{\omega \in \Omega; \ B^{i}(\omega) \
\mbox{and} \ B^{j}(\omega) \ \mbox{are continuous}\}$ ($P(\tilde
\Omega_{i,j})=1$).

Let $\omega \in \tilde \Omega_{i,j}$ and $\varepsilon>0$ fixed. For
any $(s_1,s_2) \in [0,t]^2$, define
$$\eta_{\delta}(s_1,s_2)=\left\{
\begin{array}{ll}
\alpha_H \int_0^t
\int_0^t\varphi_{\delta}(t-s_1-r_1)\varphi_{\delta}(t-s_2-r_2)|r_1-r_2|^{2H-2}dr_1
dr_2 &  \mbox{if $H>1/2$} \\
\\
\int_0^t \varphi_{\delta}(t-s_1-r)\varphi_{\delta}(t-s_2-r)dr &
\mbox{if $H=1/2$} \end{array} \right. .$$
By direct calculation, using Fubini's theorem and the fact that
$$\int_{\bR^d}
\int_{\bR^d}p_{\varepsilon}(x+B_{s_1}^i-y_1)p_{\varepsilon}(x+B_{s_2}^j-y_2)
f(y_1-y_2)dy_2dy_1=\int_{\bR^d}p_{2\varepsilon}(B_{s_1}^i-B_{s_2}^j-y)f(y)dy,$$
(which can be proved by observing that
$p_{\varepsilon_1}*p_{\varepsilon_2}=p_{\varepsilon_1+\varepsilon_2}$),
it follows that
\begin{eqnarray}
\nonumber \langle A^{\varepsilon,\delta,B^i},
A^{\varepsilon,\delta,B^j} \rangle_{\cH \cP} &=&\int_0^t \int_0^t
\int_{\bR^d}\eta_{\delta}(s_1-s_2)p_{2\varepsilon}(B_{s_1}^i-B_{s_2}^j-y)f(y)dy
ds_1 ds_2\\
\label{inner-product-A-equal-L}
&=&L_{t,2\varepsilon}^{B^i,B^j}(\eta_{\delta}).
\end{eqnarray}

\noindent Note that, for any continuous function $g:[0,t]^2 \to
\bR$,
$$\lim_{\delta \downarrow 0}\int_0^t \int_0^t
\eta_{\delta}(s_1,s_2)g(s_1,s_2)ds_1 ds_2 = \left\{
\begin{array}{ll}
\alpha_H \int_0^t \int_0^t |s_1-s_2|^{2H-2}g(s_1,s_2)ds_1 ds_2 &
\mbox{if $H>1/2$} \\
\int_0^t g(s,s)ds  & \mbox{if $H=1/2$}
\end{array} \right. $$


 In particular, we consider the (random) function
$g_{2\varepsilon}$ defined by:
$g_{2\varepsilon}(s_1,s_2)=\int_{\bR^d}p_{2\varepsilon}(B_{s_1}^i-
B_{s_2}^j-y)f(y)dy=(p_{2\varepsilon} * f)(B_{s_1}^i- B_{s_2}^j)$,
for $(s_1,s_2) \in [0,t]^2$. (Note that $g_{2\varepsilon}$ is
continuous by Remark \ref{cont-of-g}.)
Then, $$\langle A^{\varepsilon,\delta,B^i},
A^{\varepsilon,\delta,B^j} \rangle_{\cH \cP}=\int_0^t \int_0^t
\eta_{\delta}(s_1,s_2)g_{2\varepsilon}(s_1,s_2)ds_1 ds_2,$$ and
\begin{eqnarray*}
\lim_{\delta \downarrow 0} \langle A^{\varepsilon,\delta,B^i},
A^{\varepsilon,\delta,B^j} \rangle_{\cH \cP}&=&\left\{
\begin{array}{ll}
\alpha_H \int_0^t \int_0^t
|s_1-s_2|^{2H-2}g_{2\varepsilon}(s_1,s_2)ds_1 ds_2 &
\mbox{if $H>1/2$} \\
\int_0^t g_{2\varepsilon}(s,s)ds  & \mbox{if $H=1/2$}
\end{array} \right. \\
& = & L_{t,2\varepsilon}^{B^i,B^j}.
\end{eqnarray*}

{\em Step 3.} We prove that for any $(t,x) \in \bR_{+} \times \bR^d$
with $t<t_0(2)$,
\begin{equation}
\label{unif-integr-A-epsilon-delta} \{\exp(\langle
A^{\varepsilon,\delta,B^i}, A^{\varepsilon,\delta,B^j} \rangle_{\cH
\cP}) \}_{\varepsilon,\delta>0} \quad \mbox{is uniformly
integrable},
\end{equation}

Suppose first that $H=1/2$. Then $\eta_{\delta}$ satisfies condition
(\ref{cond-C1}) with $\|\eta_{\delta} \|_{1,t} \leq 1$ (see p. 318
of \cite{hu-nualart08}). By applying Proposition \ref{exp-L-bound-2}
and using (\ref{inner-product-A-equal-L}), it follows that  if $f$
is the Riesz or Bessel kernel,
$$\sup_{\varepsilon>0}E\left[\exp\left(\lambda \langle
A^{\varepsilon,\delta,B^i}, A^{\varepsilon,\delta,B^j} \rangle_{\cH
\cP}\right) \right] \leq C_{\alpha,d}^* \Phi\left(\lambda
D(t),1-\frac{d-\alpha}{2}\right), \quad \forall \lambda>0,$$ whereas
if $f$ is the heat or Poisson kernel,
$$\sup_{\varepsilon>0}E\left[\exp\left(\lambda \langle
A^{\varepsilon,\delta,B^i}, A^{\varepsilon,\delta,B^j} \rangle_{\cH
\cP}\right) \right] \leq e^{\lambda C(t)}, \quad \forall
\lambda>0.$$

\noindent Note that both constants $D(t)$ and $C(t)$ depend
(linearly) on $\|\eta_{\delta}\|_{1,t}$ (which is bounded by $1$),
and the function $\Phi(x,a)$ is increasing in $x$. We infer that
there exists an upper bound for the above supremum over
$\varepsilon$, which does not depend on $\delta$. More precisely,
denoting by $D(t),C(t)$ the respective constants $D(t),C(t)$, in
which $\|\eta_{\delta}\|_{1,t}$ is replaced by $1$, we infer that if
$f$ is the Riesz or the Bessel kernel,
\begin{equation}
\label{sup-exp-A-epsilon-delta-RB-1/2}
\sup_{\varepsilon,\delta>0}E\left[\exp\left(\lambda \langle
A^{\varepsilon,\delta,B^i}, A^{\varepsilon,\delta,B^j} \rangle_{\cH
\cP}\right) \right] \leq C_{\alpha,d}^* \Phi\left(\lambda
D(t),1-\frac{d-\alpha}{2}\right), \quad \forall \lambda>0,
\end{equation}
whereas if $f$ is the heat or the Poisson kernel,
\begin{equation}
\label{sup-exp-A-epsilon-delta-HP-1/2}
\sup_{\varepsilon,\delta>0}E\left[\exp\left(\lambda \langle
A^{\varepsilon,\delta,B^i}, A^{\varepsilon,\delta,B^j} \rangle_{\cH
\cP}\right) \right] \leq e^{\lambda C(t)}, \quad \forall \lambda>0.
\end{equation}

\noindent Relation (\ref{unif-integr-A-epsilon-delta}) follows by a
well-known criterion (see p. 218 of \cite{billingsley95}) by taking
an arbitrary $\lambda>1$.

Suppose now that $H>1/2$. Then $\eta_{\delta}$ satisfies condition
(\ref{cond-C2}) (see relation (5.13) of \cite{hu-nualart08}). By
applying Proposition \ref{exp-L-bound} and using
(\ref{inner-product-A-equal-L}), it follows that, if $f$ is the
Riesz or the Bessel kernel,
$$
\sup_{\varepsilon>0}E\left[\exp\left(\lambda \langle
A^{\varepsilon,\delta,B^i}, A^{\varepsilon,\delta,B^j} \rangle_{\cH
\cP}\right) \right] \leq C_{H,d,\alpha}^* \Phi\left(\lambda
D(t),1-\frac{d-\alpha}{2}\right),$$ for all
$0<\lambda<\lambda_0(t)$, whereas if $f$ is the heat or the Poisson
kernel,
$$
\sup_{\varepsilon>0}E\left[\exp\left(\lambda \langle
A^{\varepsilon,\delta,B^i}, A^{\varepsilon,\delta,B^j} \rangle_{\cH
\cP}\right) \right] \leq e^{\lambda C(t)}, \quad \forall
\lambda>0.$$

\noindent The constants $D(t)$ and $C(t)$ depend on $\gamma$, which
depends only on $H$. From here, we infer that, if $f$ is the Riesz
kernel or the Bessel kernel,
\begin{equation}
\label{sup-exp-A-epsilon-delta-RB}
\sup_{\varepsilon,\delta>0}E\left[\exp\left(\lambda \langle
A^{\varepsilon,\delta,B^i}, A^{\varepsilon,\delta,B^j} \rangle_{\cH
\cP}\right) \right] \leq C_{H,d,\alpha}^* \Phi\left(\lambda D(t),
1-\frac{d-\alpha}{2}\right)<\infty,
\end{equation}
for all $0<\lambda<\lambda_0(t)$, whereas if $f$ is the heat kernel
of the Poisson kernel,
\begin{equation}
\label{sup-exp-A-epsilon-delta-HP}
 \sup_{\varepsilon,\delta>0}E\left[\exp\left(\lambda \langle
A^{\varepsilon,\delta,B^i}, A^{\varepsilon,\delta,B^j} \rangle_{\cH
\cP}\right) \right] \leq e^{\lambda C(t)}<\infty, \quad \forall
\lambda>0.
\end{equation}

\noindent Relation (\ref{unif-integr-A-epsilon-delta}) follows as
before, noting that $1<\lambda_0(t)$ (since $t<t_0(2)$).

Note that (\ref{sup-u-epsilon-delta-power-k}) is obtained as a
by-product of (\ref{moment-k-u-epsilon-delta}) and
(\ref{sup-exp-A-epsilon-delta-RB-1/2})-(\ref{sup-exp-A-epsilon-delta-HP}),
since:
\begin{equation}
\label{bound-for-sup-u}
\sup_{\varepsilon,\delta>0}E[(u_{t,x}^{\varepsilon,\delta})^k] \leq
\|u_0\|_{\infty}^k \sup_{\varepsilon,\delta>0}
E\left[\exp\left(\frac{k(k-1)}{2} \langle
A^{\varepsilon,\delta,B^1}, A^{\varepsilon,\delta,B^2} \rangle_{\cH
\cP}\right) \right],
\end{equation}
 and
$$\frac{k(k-1)}{2} <\lambda_0(t) \quad \mbox{if and only if} \
\quad t<t_0(k).$$

Moreover, using
(\ref{sup-exp-A-epsilon-delta-RB-1/2})-(\ref{bound-for-sup-u}), and
the fact that $D(t),C(t)$ are increasing functions of $t$ and
$\lambda_0(t)$ is a decreasing function of $t$, we conclude that,
for any $0<T<t_0(k)$ and for any $(t,x) \in [0,T] \times \bR^d$,
$$\sup_{\varepsilon,\delta>0}E[(u_{t,x}^{\varepsilon,\delta})^k] \leq
\|u_0\|_{\infty}^k  C_{H,d,\alpha}^* \Phi\left(\frac{k(k-1)}{2}
D(T), 1-\frac{d-\alpha}{2}\right),$$
 if $f$ is the Riesz kernel or the Bessel kernel, and
$$\sup_{(t,x) \in [0,T] \times \bR^d}\sup_{\varepsilon,\delta>0}
E[(u_{t,x}^{\varepsilon,\delta})^k] \leq \|u_0\|_{\infty}^k
\exp\left\{\frac{k(k-1)}{2}C(T)\right\},$$ if $f$ is the heat or the
Bessel kernel. Hence, for any $0<T<t_0(k)$,
\begin{equation}
\label{unif-bound-u} \sup_{(t,x) \in [0,T] \times
\bR^d}\sup_{\varepsilon,\delta>0}
E[(u_{t,x}^{\varepsilon,\delta})^k] \leq K_T(k)<\infty,
\end{equation}
where $K_T(k)$ is a constant depending on $u_0,k,H,d,\alpha$ and
$T$.

{\em Step 4.} We prove that for any $0<t<t_0(k)$,
\begin{equation}
\label{limit-u-ed-power-k} \lim_{\varepsilon \downarrow 0}
\lim_{\delta \downarrow 0}
E[(u_{t,x}^{\varepsilon,\delta})^k]=\gamma_k(t,x).
\end{equation}

\noindent First we note that, for any $\varepsilon>0$
$$\lim_{\delta \downarrow 0}
E[(u_{t,x}^{\varepsilon,\delta})^k]=E\left[\prod_{j=1}^{k}u_0(x+B_t^j)\exp
\left(\sum_{1 \leq i<j \leq k}
L_{t,2\varepsilon}^{B^i,B^j}\right)\right]:=H_{t,x,k}(\varepsilon).$$

\noindent (This follows by Theorem 16.14 of \cite{billingsley95},
using (\ref{moment-k-u-epsilon-delta}), (\ref{limit-only-delta}),
and (\ref{unif-integr-A-epsilon-delta})). Next, we show that:
$$\lim_{\varepsilon \downarrow 0}
H_{t,x,k}(\varepsilon)=\gamma_{k}(t,x).$$ For this, let
$(\varepsilon_n) \downarrow 0$ be an arbitrary sequence. We have to
show that:
\begin{equation}
\label{lim-mom-k-only-epsilon} \lim_{n \to \infty}
H_{t,x,k}(\varepsilon_n)=\gamma_{k}(t,x).
\end{equation}

We use the fact that in a metric space, $x_n \to x$ if and only if
for any subsequence $N' \subset \bN $ there exists a sub-subsequence
$N'' \subset N'$ such that $x_{n} \to x$, as $n \to \infty, n \in
N''$. Let $N' \subset \bN$ be an arbitrary subsequence.

By Lemma \ref{existence-of-L}, the limit
$L_{t}^{B^i,B^j}:=\lim_{\varepsilon \downarrow
0}L_{t,2\varepsilon}^{B^i,B^j}$ exists in $L^2(\Omega)$. (In this
lemma, we take $\eta(r,s)=\alpha_H|r-s|^{2H-2}$ if $H>1/2$, and
$\eta(r,s)=1_{\{r=s\}}$ if $H=1/2$.) Hence
$L_{t,2\varepsilon_{n}}^{B^i,B^j} \to L_{t}^{B^i,B^j}$ in
probability, as $n \to \infty,  n \in N'$, and there exists a
sub-subsequence $N'' \subset N'$ such that
$$L_{t,2\varepsilon_{n}}^{B^i,B^j} \to L_{t}^{B^i,B^j} \quad \mbox{a.s.} \quad
\mbox{as}  \ n \to \infty,n \in N''.$$

\noindent Note that $E[\exp (\lambda
L_{t,2\varepsilon}^{B^i,B^i})]<\infty$ for all $\lambda>0$,
respectively for all $0<\lambda<\lambda_0(t)$, with
$\lambda_0(t)>1$. (If $H=1/2$, we use Proposition
\ref{exp-L-bound-2} with $\eta(r,s)=1_{\{r=s\}}$. If $H>1/2$, we use
Proposition \ref{exp-L-bound} with
$\eta(r,s)=\alpha_H|r-s|^{2H-2}$.) Hence,
$\{\exp(L_{t,2\varepsilon}^{B^{i},B^j}) \}_{\varepsilon
>0}$ is uniform integrable. By Theorem 16.14 of \cite{billingsley95},
$$H_{t,x,k}(\epsilon_n) \to \gamma_k(t,x), \quad \mbox{as} \quad
n \to \infty, n \in N''.$$ Relation (\ref{lim-mom-k-only-epsilon})
follows using the above-mentioned subsequence criterion.

{\em Step 5.} We prove that for any $0<t<t_0(2)$,
\begin{equation}
\label{limit-u-ed-prime} \lim_{\varepsilon,\varepsilon' \downarrow
0} \lim_{\delta, \delta' \downarrow
0}E[u_{t,x}^{\varepsilon,\delta}u_{t,x}^{\varepsilon',\delta'}]=\gamma_2(t,x).
\end{equation}

Similarly to (\ref{moment-k-u-epsilon-delta}) and
(\ref{limit-only-delta}), one can prove that:
$$E[u_{t,x}^{\varepsilon,\delta}u_{t,x}^{\varepsilon',\delta'}]=
E[u_0(x+B_t^1)u_0(x+B_t^2) \exp(\langle A^{\varepsilon,\delta, B^1},
A^{\varepsilon',\delta',B^2} \rangle_{\cH \cP})],$$ and
$$\lim_{\delta,\delta'
\downarrow 0}\langle A^{\varepsilon,\delta, B^1},
A^{\varepsilon',\delta',B^2} \rangle_{\cH
\cP}=L_{t,\varepsilon+\varepsilon'}^{B^1,B^2} \ \ \forall
\varepsilon>0, \ \forall \varepsilon'>0 \ \mbox{a.s.}$$

\noindent Relation (\ref{limit-u-ed-prime}) follows using the same
argument as in {\em Step 4} (based on
(\ref{unif-integr-A-epsilon-delta})).

{\em Step 6.} We prove that for any $t<t_0(2)$ and $x \in \bR^d$,
the limit
$$u_{t,x}:=\lim_{\varepsilon
\downarrow 0} \lim_{\delta \downarrow 0}u_{t,x}^{\varepsilon,\delta}
\ \mbox{exists in} \ L^2(\Omega).$$

Let $t<t_0(2)$ and $x \in \bR^d$ be fixed. From
(\ref{limit-u-ed-prime}), we obtain that:
\begin{equation}
\label{limit-u-ed-prime-2} \lim_{\varepsilon,\varepsilon' \downarrow
0} \lim_{\delta, \delta' \downarrow
0}E|u_{t,x}^{\varepsilon,\delta}-u_{t,x}^{\varepsilon',\delta'}|^2=0.
\end{equation}

Let $(\varepsilon_n)_n \downarrow 0$ and $(\delta_n)_n \downarrow 0$
be arbitrary. From (\ref{limit-u-ed-prime-2}), it follows that
$(u_{t,x}^{\varepsilon_n,\delta_n})_n$ is a Cauchy sequence in
$L^2(\Omega)$. Hence, there exists $u_{t,x} \in L^2(\Omega)$ such
that $E|u_{t,x}^{\varepsilon_n,\delta_n}-u_{t,x}|^2 \to 0$. The fact
that $u_{t,x}$ does not depend on $(\varepsilon_n)_n$ and
$(\delta_n)$ is proved by a standard argument (see for instance, the
proof of Lemma \ref{existence-of-L}).

{\em Step 7.} We now prove (\ref{u-power-k}).

If $t<t_0(M)$ for some $M \geq 3$, then
$\sup_{\varepsilon,\delta>0}E[(u_{t,x}^{\varepsilon,\delta})^M]<\infty$,
and hence,
$\{(u_{t,x}^{\varepsilon,\delta})^p\}_{\varepsilon,\delta>0}$ is
uniformly integrable, for any $2 \leq p<M$. Since
$u_{t,x}^{\varepsilon,\delta}-u_{t,x} \to 0$ a.s. (along a
subsequence), we conclude that
$E|u_{t,x}^{\varepsilon,\delta}-u_{t,x}|^p \to 0$. In particular,
$\lim_{\varepsilon \downarrow 0} \lim_{\delta \downarrow
0}E[(u_{t,x}^{\varepsilon,\delta})^k] = E[u_{t,x}^k]$ for any $2
\leq k<M$, and (\ref{u-power-k}) follows by
(\ref{limit-u-ed-power-k}).

{\em Step 8.} We prove that $u=\{u_{t,x}; (t,x) \in [0,t_0(2))
\times \bR^d\}$ is a solution of (\ref{heat}).

Let $(u_{t,x}^{\varepsilon,\delta})_{t,x}$ be the solution of
equation (\ref{heat-epsilon-delta}) (whose existence is guaranteed
by Proposition \ref{exist-sol-ed}). Let $(\varepsilon_n) \downarrow
0$ and $(\delta_n) \downarrow 0$ be arbitrary. By {\em Step 6},
$u_{t,x}^{\varepsilon_n,\delta_n} \to u_{t,x}$ in $L^2(\Omega)$, and
hence $u_{t,x}^{\varepsilon_n,\delta_n} \to u_{t,x}$ a.s. along a
subsequence $N'\subset \bN$. Since
$u_{t,x}^{\varepsilon_n,\delta_n}$ is $\cF_t$-measurable, it follows
that $u_{t,x}$ is $\cF_t$-measurable.

Let $(t,x) \in \bR_{+} \times \bR^d$ be fixed with $t<t_0(2)$. We
have to show that (\ref{def-sol-u}) holds. As in the proof of
Proposition \ref{exist-sol-ed}, it suffices to show that
(\ref{def-sol-u}) holds for $F=F_{\varphi}$ with $\varphi \in \cH
\cP$. Moreover, it suffices to take $\varphi=1_{[0,a]}\phi$, with $a
\in \bR_{+}, \phi \in \cP(\bR^d)$, since the class of these
functions is dense in $\cH \cP$.

The idea is to take the limit in (\ref{def-sol-u-e-d}), as
$\varepsilon \downarrow 0, \delta \downarrow 0$. On the left-hand
side, $E(F_{\varphi} u_{t,x}^{\varepsilon,\delta}) \to E(F_{\varphi}
u_{t,x})$, since $E|u_{t,x}^{\varepsilon,\delta}-u_{t,x}|^2 \to 0$.
On the right-hand side of (\ref{def-sol-u-e-d}), we have: (see
({\ref{FK-formula}))
\begin{eqnarray*}
E \langle Y^{t,x,\varepsilon,\delta}, DF_{\varphi} \rangle_{\cH \cP}
&=& \int_0^{t} \int_{\bR^d} p_{t-s}(x-y)
E(u_{s,y}^{\varepsilon,\delta}F_{\varphi})
V^{\varepsilon,\delta}(s,y)dyds.
\end{eqnarray*}

We have $\lim_{\varepsilon \downarrow 0}\lim_{\delta \downarrow 0}
E(F u_{s,y}^{\varepsilon,\delta})=E(F_{\varphi} u_{s,y})$, for all
$(s,y) \in [0,t] \times \bR^d$. Note that
$V^{\varepsilon,\delta}=\psi_{\varepsilon,\delta}*H$, where
$\psi_{\varepsilon,\delta}(s,y)=\varphi_{\delta}(s)
p_{\varepsilon}(y)$ and
$$H(s,y):=\alpha_{H} \int_{\bR_+} \int_{\bR^d}
\varphi(s',y')|s-s'|^{2H_2}f(y-y')dy'ds'.$$

\noindent Hence, $\lim_{\varepsilon \downarrow 0} \lim_{\delta
\downarrow 0} V^{\varepsilon,\delta}(s,y)= H(s,y)$, for all $(s,y)
\in [0,t] \times \bR^d$. Therefore, by the dominated convergence
theorem (whose application is justified below),
$$\lim_{\varepsilon \downarrow 0} \lim_{\delta \downarrow
0} E \langle Y^{t,x,\varepsilon_n,\delta_n}, DF_{\varphi}
\rangle_{\cH \cP}= \int_0^t
\int_{\bR^d}p_{t-s}(x-y)E(u_{s,y}F_{\varphi}) H(s,y)dyds.$$ A direct
calculation shows that the limit above coincides with $E \langle
Y^{t,x}, DF_{\varphi} \rangle_{\cH \cP}$ (using the fact that
$D_{s,y}F_{\varphi}=\varphi(s,y)F_{\varphi}$).

It remains to justify the application of the dominated convergence
theorem. Using the Cauchy-Schwartz inequality and
(\ref{unif-bound-u}), we have:
\begin{equation}
\label{bound-for-uF} \sup_{(s,y) \in [0,t] \times \bR^d}
\sup_{\varepsilon,\delta>0}
|E(u_{s,y}^{\varepsilon,\delta}F_{\varphi})|  \leq
\{K_t(2)\}^{1/2}\{E(F_{\varphi}^2)\}^{1/2}:=K_t^*.
\end{equation}

Note that
$V^{\varepsilon,\delta}(s,y)=J^{\delta}(s)I_{\varepsilon}(y)$, where
$$J_{\delta}(s)=\langle \varphi_{\delta}(s- \cdot),1_{[0,a]}
\rangle_{\cH} \quad \mbox{and} \quad I_{\varepsilon}(y) =\langle
p_{\varepsilon}(y-\cdot),\phi \rangle_{\cP(\bR^d)}.$$

\noindent Let $\beta=(\beta_t)_{t \geq 0}$ be a fBm of index $H$.
Then
\begin{eqnarray*}
J_{\delta}(s) &=& =\frac{1}{\delta} \alpha_{H} \int_{s-\delta}^{s}
\int_0^a |r-r'|^{2H-2}dr dr' =\frac{1}{\delta}
E[(\beta_s-\beta_{s-\delta})\beta_{a}] \\
&=& \frac{1}{2\delta}
[s^{2H}+|s-\delta-a|^{2H}-|s-a|^{2H}-|s-\delta|^{2H}].
\end{eqnarray*}

\noindent Using an argument similar to (5.13) of
\cite{hu-nualart08}, one can show that:
\begin{equation}
\label{bound-for-J} J_{\delta}(s) \leq 2H s^{2H-1}+ c_{a,H}, \quad
\forall s \in [0,t], \forall \delta >0,
\end{equation}
where $c_{a,H}=H (a^{2H-1}+1)$. (This argument is based on treating
separately the cases $a>s,a \leq s$, and considering in each case
several intervals for $\delta$.)

We claim that:
\begin{equation}
\label{bound-for-I} |I_{\varepsilon}(y)| \leq c_{\phi}, \quad
\forall y \in \bR^d, \forall \varepsilon>0,
\end{equation}
where $c_{\phi}$ is a constant depending on $d,f$ and $\phi$. To see
this, we assume without loss of generality that $\phi=\psi*p_{b}$,
for some $\psi \in C_0^{\infty}(\bR^d),b>0$ (since functions of this
form are dense in $\cP(\bR^d)$). Then, assuming that
$\mu(d\xi)=g(\xi)d\xi$, we have:
\begin{eqnarray*}
|I_{\varepsilon}(y)| &=& \left|\int_{\bR^d} \overline{\cF
p_{\varepsilon}(y-\cdot)(\xi)} \cF \phi (\xi)g(\xi)d\xi \right| \leq
\int_{\bR^d} |\cF p_{\varepsilon}(y-\cdot)(\xi)| \ |\cF \phi (\xi)|
g(\xi)d\xi \\
&=& (2\pi)^{-2d} \int_{\bR^d} e^{-\varepsilon |\xi|^2/2} |\cF
\psi(\xi)|e^{-b|\xi|^2/2} g(\xi)d\xi \\
& \leq & (2\pi)^{-2d} \int_{\bR^d}  |\cF \psi(\xi)|e^{-b|\xi|^2/2}
g(\xi)d\xi \\
& \leq & (2\pi)^{-2d} \left(\int_{\bR^d} |\cF
\psi(\xi)|^2g(\xi)d\xi\right)^{1/2}
\left(\int_{\bR^d} e^{-b|\xi|^2}g(\xi)d\xi \right)^{1/2} \\
&=& (2\pi)^{-3d/2} \|\psi\|_{\cP(\bR^d)} \
[\psi^{(1)}(b,b)]^{1/2}:=c_{\phi},
\end{eqnarray*}
where we used the Cauchy-Schwartz inequality for the last inequality
above.

Combining (\ref{bound-for-J}) and (\ref{bound-for-I}), we get:
\begin{equation}
\label{bound-for-V} \sup_{(s,y) \in [0,t] \times \bR^d}
\sup_{\varepsilon,\delta>0}|V^{\varepsilon,\delta}(s,y)| \leq
c_{\phi}(2Hs^{2H-1}+c_{a,H}).
\end{equation}

\noindent From (\ref{bound-for-uF}) and (\ref{bound-for-V}), we
infer that for any $\varepsilon>0,\delta>0,(s,y) \in [0,t] \times
\bR^d$,
$$|p_{t-s}(x-y) E(u_{s,y}^{\varepsilon,\delta}F_{\varphi})
V^{\varepsilon,\delta}(s,y)| \leq \Psi(s,y),$$ where
$\Psi(s,y):=K_{t}^{*}c_{\phi} p_{t-s}(x-y)(2H s^{2H-1} +c_{a,H})$.
Finally, we note that $\Psi$ is integrable on $[0,t] \times \bR^d$,
since:
\begin{eqnarray*}
\int_{0}^{t} \int_{\bR^d}\Psi(s,y)dy ds &=& K_t^*c_{\phi} \int_0^t
(2H
s^{2H-1} +c_{a,H}) \left( \int_{\bR^d}p_{t-s}(x-y)dy \right)ds \\
&=&  K_t^*c_{\phi} (t^{2H}+c_{a,H}t)<\infty.
\end{eqnarray*}
$\Box$

 \noindent \small{{\bf Acknowledgement.} The
authors are grateful to an anonymous referee who read the article
very carefully and made numerous suggestions for improving the
presentation and clarifying some technical points.}

\normalsize{

\appendix

\section{Correction to Theorem 3.13 of \cite{balan-tudor08}}

Theorem 3.13 of \cite{balan-tudor08} gives the necessary and
sufficient condition for the existence of the solution of the
equation: $u_t=\frac{1}{2}\Delta u+\dot W$ in $(0,T) \times \bR^d$,
with $u(0,\cdot)=0$. This condition is equivalent to saying that
$\|g_{tx}\|_{\cH \cP}<\infty$, where $g_{tx}(s,y)=[2\pi
(t-s)]^{-d/2} \exp\{-|x-y|^2/[2(t-s)]\}=p_{t-s}(x-y)$. The condition
is incorrectly stated in the case of the Bessel kernel, the heat
kernel, and the Poisson kernel. We state below the correction of
this result, whose proof will appear as an erratum in
\cite{balan-tudor08-err}. In connection with the present article, we
observe that:
$$\|g_{tx}\|_{\cH \cP}^2 = \alpha_H \int_0^t \int_0^t
|r-s|^{2H-2}I_{f}(r,s)drds=\alpha_1(t),$$ since
$I_{f}(r,s):=\int_{\bR^d} \int_{\bR^d}g_{tx}(s,y)f(y-z)g_{tx}(r,z)dy
dz=\psi^{*(1)}(r,s)$.

\begin{theorem}
(i) If $f$ is the Riesz kernel of order $\alpha$, or the Bessel
kernel of order $\alpha$, then $\|g_{tx}\|_{\cH \cP}<\infty$ if and
only if $H>(d-\alpha)/4$.

(ii) If $f$ is the heat kernel of order $\alpha$, or the Poisson
kernel of order $\alpha$, then $\|g_{tx}\|_{\cH \cP}<\infty$ for any
$H>1/2$ and $d \geq 1$.
\end{theorem}

}


\begin{thebibliography}{99}


\bibitem{balan-tudor08} Balan, R. M. and Tudor, C. A. (2008). The
stochastic heat equation with fractional-colored noise: existence of
the solution. {\em Latin Amer. J. Probab. Math. Stat.} {\bf 4},
57-87.

\bibitem{balan-tudor08-err} Balan, R. M. and Tudor, C. A. (2009). Erratum to:
``The stochastic heat equation with fractional-colored noise:
existence of the solution''. {\em Latin Amer. J. Probab. Math.
Stat.}

\bibitem{BPS04} Bayraktar, E., Poor, V. and Sircar, R. (2004).
Estimating the fractal dimension of the SP 500 index using wavelet
analysis. {\em  Intern. J. Theor. Appl. Finance} {\bf 7}, 615-643.

\bibitem{billingsley68} Billingsley, P. (1968). {\em Convergence of Probability
Measures}. John Wiley, New York.


\bibitem{billingsley95} Billingsley, P. (1995). {\em Probability and
Measure}. Third edition. John Wiley, New York.

\bibitem{buckdahn-nualart94} Buckdahn, R. and Nualart, D.
(1994). Linear stochastic differential equations and Wick products.
{\em Probab. Th. Rel. Fields} {\bf 99}, 501-526.

\bibitem{CCL03} del Castillo-Negrete. D, Carreras, B. A. and Lynch,
V. E. (2003). Front dynamics in reaction-diffusion systems with Levy
flights: a fractional difussion approach. {\em Phys. Rev. Letters}
{\bf 91}.

\bibitem{CNT} Coutin, L.,  Nualart, D.,  Tudor, C.A. (2001). The Tanaka formula
for the fractional Brownia motion.  {\em Stochastic Proc. Applic.} {
\bf 94}, 301-315.

\bibitem{dalang99} Dalang, R. C. (1999). Extending martingale measure stochastic
integral with application to spatially homogenous s.p.d.e.'s. {\em
Electr. J. Probab.} {\bf 4}, paper 6, 1-29.

\bibitem{DaFr98} Dalang, R. C. and Frangos, N. (1998). The stochastic
wave equation in two spatial dimensions. {\em Ann. Probab.} {\bf
26}, 187-212.

\bibitem{DaMu03} Dalang, R. C. and Mueller, C. (2003). Some non-linear
S.P.D.E.'s that are second order in time. {\em Electr. J. Probab.}
{\bf 8}, paper 1, 1-21.


\bibitem{DaSa05} Dalang, R. C. and Sanz-Sol\'e, M. (2005). Regularity of
the sample paths of a class of second order s.p.d.e.'s. {\em J.
Funct. Anal.} {\bf 227}, 304-337.

\bibitem{DMS03} Denk, G., Meintrup, D. and Schaffer, S. (2004).
Modeling, simulation and optimization of integrated circuits. {\em
Intern. Ser. Numerical Math.} {\bf 146}, 251-267.

\bibitem{DMD}
Duncan, T. E., Maslowski, B. and  Pasik-Duncan, B. (2002).
Fractional Brownian motion and stochastic equations in Hilbert
spaces. {\em Stoch. Dyn.} {\bf 2}, 225-250.

\bibitem{folland95} Folland, G.B. (1995). {\em Introduction to
Partial Differential Equations}, Second Edition.  Princeton
University Press, Princeton.


\bibitem{KFN08} Foondun, M., Khosnevisan, D. and Nualart, E. (2008).
A local-time correspondence for stochastic partial differential
equations. Preprint.

\bibitem{hu01} Hu, Y. (2001). Heat equations with fractional white
noise potentials. {\em Appl. Math. Optim.} {\bf 43}, 221-243.

\bibitem{hu-nualart08} Hu, Y. and Nualart, D. (2009). Stochastic
heat equation driven by fractional noise and local time. {\em
Probab. Theory Rel. Fields} {\bf 143}, 285-328.


\bibitem{karatzas-shreve91} Karatzas, I. and Shreve, S.E. (1991). {\em
Brownian Motion and Stochastic Calculus}, Second edition, Springer,
New York.


\bibitem{KouSinney04} Kou, S.C. and Sunney, X. (2004). Generalized
Langevin equation with fractional Gaussian noise: subdiffusion
within a single protein molecule. {\em Phys. Rev. Letters} {\bf
93}(18).

\bibitem{leon-sanmartin07} Le\'{o}n, J. A. and San Mart\'{i}n,
J. (2007). Linear stochastic differential equations driven by a
fractional Brownian motion with Hurst parameter less than $1/2$.
{\em Stoch. Anal. Appl.} {\bf 25}, 105-126.


\bibitem{MMV01} M\'emin, J., Mishura, Y. and Valkeila, E. (2001).
Inequalities for the moments of Wiener integrals with respect to
fractional Brownian motions. {\em Stat. Probab. Letters} {\bf 51},
197-206.

\bibitem{MiSa99} Millet, A. and Sanz-Sol\'e, M (1999). A stochastic wave
equation in two space dimension: smoothness of the law. {\em Ann.
Probab.} {\bf 27}, 803-844.

\bibitem{MaNu}
Maslovski, B. and Nualart, D. (2003). Evolution equations driven by
a fractional Brownian motion. {\em J. Funct. Anal.} {\bf 202},
277-305.

\bibitem{nourdin-tudor06} Nourdin, I. and Tudor, C. A. (2006).
Some linear fractional stochastic equations. {\em Stochastics} {\bf
78}, 51-65.

\bibitem{nualart98} Nualart, D. (1998). Analysis on Wiener space and
anticipative stochastic calculus. {\em Lect. Notes Math.} {\bf
1690}, Springer-Verlag, Berlin.

\bibitem{nualart03} Nualart, D. (2003). Stochastic integration with respect to
fractional Brownian motion and applications. {\em Contem. Math.}
{\bf 336}, 3-39.

\bibitem{nualart06} Nualart, D. (2006). {\em Malliavin Calculus
and Related Topics}, Second Edition. Springer-Verlag, Berlin.

\bibitem{NuOuk}
Nualart, D. and Ouknine, Y. (2004). Regularization of quasilinear
heat equation equations by a fractional noise. {\em Stoch. Dyn.}
{\bf 4}, 201-221.

\bibitem{nualart-rozovskii97} Nualart, D. and Rozovskii, B.
(1997). Weighted stochastic Sobolev spaces and bilinear SPDEs driven
by space-time white noise. {\em J. Funct. Anal.} {\bf 149}, 200-225.


\bibitem{NV}
Nualart, D. and Vives, J. (1992).  Smoothness of Brownian local
times and related functionals. {\em Potential Analysis}  {\bf 1},
257-263.

\bibitem{nualart-zakai89} Nualart, D. and Zakai, M. (1989).
Generalized Brownian functionals and the solution to a stochastic
partial differential equation. {\em J. Funct. Anal.} {\bf 84},
279-296.

\bibitem{PeZa00} Peszat, S. and Zabczyk, J. (2000). Nonlinear stochastic
wave and heat equations. {\em Probab. Th. Rel. Fields} {\bf 116},
421-443.

\bibitem{TTV}
Tindel, S., Tudor, C. A. and Viens, F. (2003). Stochastic evolution
equations  with fractional Brownian motion. {\em Probab. Th. Rel.
Fields} {\bf 127}, 186-204.

\bibitem{tudor04} Tudor, C. (2004). Fractional bilinear stochastic
equations with the drift in the first fractional chaos. {\em Stoch.
Anal. Appl.} {\bf 22}, 1209-1233.

\bibitem{walsh86} Walsh, J. B. (1986). An introduction to stochastic
partial differential equations. {\em Ecole d'Et\'{e} de
Probabilit\'{e}s de Saint-Flour XIV. Lecture Notes in Math.} {\bf
1180}, 265-439. Springer-Verlag, Berlin.
\end{thebibliography}
\end{document}